\documentclass[12pt]{article} 
\usepackage[english]{babel}
\usepackage[cp1251]{inputenc}
\usepackage{amssymb,amscd,amsmath}

\usepackage[linktocpage=true, colorlinks=true, linkcolor=blue, citecolor=blue, urlcolor=blue]{hyperref}

\begin{document}
\thispagestyle{empty}

\parbox{0.275\textwidth }{
\begin{picture}(130,120)
\put(67,67){{\Huge $\frac{dx}{dt}$}}
\put(40,0){\vector(0,1){70}}
\put(40,70){\line(0,1){10}}
\qbezier(40,80)(40,110)(80,110)
\qbezier(80,110)(120,110)(120,80)
\put(120,80){\line(0,-1){10}}
\qbezier(120,70)(120,35)(75,35)
\put(75,35){\vector(-1,0){5}}
\put(0,35){\line(1,0){70}}
\put(0,40){\vector(1,0){10}}
\qbezier(10,40)(35,40)(35,70)
\put(35,70){\line(0,1){10}}
\qbezier(35,80)(35,115)(80,115)
\qbezier(80,115)(125,115)(125,80)
\put(125,80){\line(0,-1){10}}
\qbezier(125,70)(125,30)(75,30)
\put(75,30){\line(-1,0){15}}
\qbezier(60,30)(45,30)(45,15)
\put(45,15){\vector(0,-1){15}}
\end{picture}
}
\hfill
\noindent \parbox{0.45\textwidth }{\footnotesize \it
\begin{center}
DIFFERENTIAL EQUATIONS\\ 
AND\\
CONTROL PROCESSES

\noindent N. 2, 2020\\

\noindent Electronic Journal, \\
reg. N ${\Phi}$C77-39410 at 15.04.2010

\noindent ISSN 1817-2172
\medskip

\noindent http://diffjournal.spbu.ru/\\

\noindent e-mail: jodiff@mail.ru

\end{center}}
\vspace*{10mm}
\large

\rightline{\underline{\it Stochastic differential equations}}
\rightline{\underline{\it Numerical methods}}
\rightline{\underline{\it Computer modeling in dynamical and control systems}}

\bigskip

\begin{center}
\large{\bf The Proof of Convergence with Probability 1 in the Method of
Expansion of Iterated It\^{o} Stochastic Integrals 
Based on Generalized Multiple Fourier Series}
\end{center}

\centerline{\large {Dmitriy F. Kuznetsov}}
\begin{center}
Peter the Great Saint-Petersburg Polytechnic University\\
\hspace{-2mm}e-mail: sde\_\hspace{0.5mm}kuznetsov@inbox.ru
\end{center}

\vspace{5mm}

\noindent
{\bf Abstract.} The article is devoted to 
the formulation and proof of the theorem on 
convergence with probability 1 (w.~p.~1) of 
expansion of iterated It\^{o} stochastic integrals of arbitrary
multiplicity
based on generalized multiple Fourier series
converging in the sense of norm in Hilbert space.
The cases of multiple Fourier--Legendre series
and multiple trigonomertic Fourier series
are considered in detail.
The proof of the mentioned theorem is based
on the general properties of multiple Fourier
series as well as on the estimate for the fourth
moment of approximation error in the
method of expansion of iterated It\^{o} stochastic integrals
based on generalized multiple Fourier series.\\
{\bf Key words:} Iterated It\^{o} stochastic integral,
generalized multiple Fourier series,
multiple Fourier--Legendre series, multiple
trigonometric Fourier series, Parseval equality,
Legendre polynomials, convergence w.~p.~1, 
mean-square convergence, convergence in the mean of arbitrary degree,
expansion, approximation.

\newpage

\renewcommand{\baselinestretch}{1.3}

\newpage
\noindent
\section{Introduction}

The beginning of an intensive study of the problem of 
mean-square approximation of iterated It\^{o}
and Stratonovich stochastic integrals in the context of
the numerical solution of It\^{o} stochastic differential 
equations dates back to the 1980s--1990s. 
To date, there are many publications on the mentioned problem
\cite{Milstein}-\cite{Kuznetsov16} 
(also see bibliographic references in these works).
There are various approaches to solving the problem of the
mean-square
approximation of iterated stochastic integrals.
Among them, we note the approach based on the 
Karhunen--Loeve expansion of the Brownian bridge process 
\cite{Milstein}-\cite{Kloeden3}, \cite{Milstein2}, \cite{Platen},
\cite{Zahri},
approach based on the expansion of the Wiener process using 
various basis systems of functions \cite{Averina}, \cite{Prigarin},
\cite{Rybakov}, \cite{Kuznetsov14},
approach based on the conditional joint characteristic 
function of a stochastic integral of multiplicity 2 \cite{Wiktorsson1},
\cite{Wiktorsson2}
as well as an approach based on multiple integral sums
\cite{Milstein}, \cite{Allen}.

The use of multiple and iterated generalized Fourier series 
by various complete orthonormal systems of 
functions in the space $L_2([t, T])$ 
for the expansion of iterated
It\^{o} and Stratonovich stochastic integrals
was reflected in a number of author's works
\cite{Kuznetsov0}-\cite{Kuznetsov2},
\cite{Kuznetsov3}-\cite{Kuznetsov5},
\cite{Kuznetsov6},
\cite{Kuznetsov7}-\cite{Kuznetsov13a}, \cite{Kuznetsov15}.
The mentioned results based on generalized 
multiple and iterated Fourier series are systematized in the monograph
\cite{Kuznetsov16} (2022).

The idea of the method of 
expansion of iterated It\^{o} stochastic integrals
based on generalized multiple Fourier series
is as follows: the iterated It\^{o} stochastic 
integral of multiplicity $k$ $(k\in \mathbb{N})$
is represented as a multiple stochastic 
integral from the certain discontinuous nonrandom function of $k$ variables
defined on the hypercube $[t, T]^k$, where $[t, T]$ is an interval of 
integration of the iterated It\^{o} stochastic integral. Then, 
the indicated 
nonrandom function is expanded into the generalized 
multiple Fourier series converging 
in the sense of norm
in the space 
$L_2([t,T]^k)$. After a number of nontrivial transformations we come 
\cite{Kuznetsov3} (2006)
to the mean-square converging expansion of the iterated It\^{o} stochastic 
integral into the multiple 
series of products
of standard  Gaussian random 
variables. The coefficients of this 
series are the coefficients of 
generalized multiple Fourier series for the mentioned nonrandom function 
of $k$ variables, which can be calculated using the explicit formula 
regardless of multiplicity $k$ of the iterated It\^{o} stochastic integral.

In a lot of author's publications 
the convergence of the method of 
expansion of iterated It\^{o} stochastic integrals
based on generalized multiple Fourier series
has been considered in different
probabilistic
meanings. For example, the mean-square convergence 
\cite{Kuznetsov3}-\cite{Kuznetsov5},
\cite{Kuznetsov6},
\cite{Kuznetsov7}-\cite{Kuznetsov13a}, \cite{Kuznetsov15}, \cite{Kuznetsov16}
and convergence in the mean of degree
$2n$ $(n\in\mathbb{N})$ 
\cite{Kuznetsov16} (Sect.~1.1.9, 1.11, 1.12),\ \cite{arxiv-1} (Sect.~6, 15, 16)
have been proved. On the examples
of specific iterated It\^{o} stochastic integrals of mutiplicities 1 and 2
the convergence w.~p.~1 also has been considered 
\cite{Kuznetsov4}-\cite{Kuznetsov5},
\cite{Kuznetsov6},
\cite{Kuznetsov7}, \cite{Kuznetsov8}.
This article is devoted to the development of the method of 
expansion of iterated It\^{o} stochastic integrals
based on generalized multiple Fourier series. Namely,
we formulate and prove the theorem on convergence
w.~p.~1 
of the mentioned method for an arbitrary
multiplicity $k$ $(k\in\mathbb{N})$
of the iterated It\^{o} stochastic integrals.
Moreover, the cases of multiple Fourier--Legendre
series and multiple trigonometric Fourier series
are considered in detail.

\section{Method of Expansion 
of Iterated It\^{o} Stochastic Integrals of Multiplicity
$k$ $(k\in\mathbb{N})$ Based on Generalized Multiple Fourier
Series}

Let $(\Omega,$ ${\rm F},$ ${\sf P})$ be a complete probability space, let 
$\{{\rm F}_t, t\in[0,T]\}$ be a nondecreasing right-continuous 
family of $\sigma$-algebras of ${\rm F},$
and let ${\bf w}_t$ be a standard $m$-dimensional Wiener stochastic 
process, which is
${\rm F}_t$-measurable for any $t\in[0, T].$ We assume that the components
${\bf w}_{t}^{(i)}$ $(i=1,\ldots,m)$ of this process are independent.

Let us consider an efficient method 
\cite{Kuznetsov3}-\cite{Kuznetsov5},
\cite{Kuznetsov6},
\cite{Kuznetsov7}-\cite{Kuznetsov13a}, \cite{Kuznetsov15}-\cite{arxiv-1}
of the expansion and mean-square approximation 
of iterated It\^{o} stochastic integrals 
of the form
\begin{equation}
\label{sodom20}
J[\psi^{(k)}]_{T,t}
=\int\limits_t^T\psi_k(t_k) \ldots \int\limits_t^{t_{2}}
\psi_1(t_1) d{\bf w}_{t_1}^{(i_1)}\ldots
d{\bf w}_{t_k}^{(i_k)},
\end{equation}
\noindent
where $0\le t<T < \infty,$ $\psi_l(\tau)$ $(l=1,\ldots,k)$ are
nonrandom functions from the space $L_2([t, T])$,
${\bf w}_{\tau}^{(i)}$ ($i=1,\ldots,m)$ are independent
standard Wiener processes and
${\bf w}_{\tau}^{(0)}=\tau,$
$i_1,\ldots,i_k=0, 1,\ldots,m.$

Suppose that $\{\phi_j(x)\}_{j=0}^{\infty}$
is a complete orthonormal system of functions in the space
$L_2([t, T])$ and
define the following function on the hypercube $[t, T]^k$
\begin{equation}
\label{ppp}
K(t_1,\ldots,t_k)=
\begin{cases}
\psi_1(t_1)\ldots \psi_k(t_k),\ &t_1<\ldots<t_k\\
~\\
0,\ &\hbox{\rm otherwise}
\end{cases},
\end{equation}

\newpage
\noindent
where $t_1,\ldots,t_k\in [t, T]$ $(k\ge 2)$ and 
$K(t_1)\equiv\psi_1(t_1)$ for $t_1\in[t, T].$

The function $K(t_1,\ldots,t_k)$ belongs to the space 
$L_2([t, T]^k).$
At this situation it is well known that the generalized 
multiple Fourier series 
of $K(t_1,\ldots,t_k)\in L_2([t, T]^k)$ converges
to $K(t_1,\ldots,t_k)$ on the hypercube $[t, T]^k$ in 
the mean-square sense, i.e.
\begin{equation}
\label{sos1z}
\hbox{\vtop{\offinterlineskip\halign{
\hfil#\hfil\cr
{\rm lim}\cr
$\stackrel{}{{}_{p_1,\ldots,p_k\to \infty}}$\cr
}} }\left\Vert
K(t_1,\ldots,t_k)-
\sum_{j_1=0}^{p_1}\ldots \sum_{j_k=0}^{p_k}
C_{j_k\ldots j_1}\prod_{l=1}^{k} \phi_{j_l}(t_l)\right
\Vert_{L_2([t, T]^k)}=0,
\end{equation}
where
\begin{equation}
\label{ppppa}
C_{j_k\ldots j_1}=\int\limits_{[t,T]^k}
K(t_1,\ldots,t_k)\prod_{l=1}^{k}\phi_{j_l}(t_l)dt_1\ldots dt_k
\end{equation}
is the Fourier coefficient and
$$
\left\Vert f\right\Vert_{L_2([t, T]^k)}=\left(\int\limits_{[t,T]^k}
f^2(t_1,\ldots,t_k)dt_1\ldots dt_k\right)^{1/2}.
$$

Consider the discretization $\{\tau_j\}_{j=0}^N$ of $[t,T]$ such that
\begin{equation}
\label{1111}
t=\tau_0<\ldots <\tau_N=T,\ \ \ \
\Delta_N=
\hbox{\vtop{\offinterlineskip\halign{
\hfil#\hfil\cr
{\rm max}\cr
$\stackrel{}{{}_{0\le j\le N-1}}$\cr
}} }\Delta\tau_j\to 0\ \ \ \hbox{if}\ \ \ N\to \infty,
\end{equation}
where $\Delta\tau_j=\tau_{j+1}-\tau_j.$

{\bf Theorem 1} \cite{Kuznetsov3} (2006), \cite{Kuznetsov4}-\cite{Kuznetsov5},
\cite{Kuznetsov6},
\cite{Kuznetsov7}-\cite{Kuznetsov13a}, \cite{Kuznetsov15}, 
\cite{Kuznetsov16}.
{\it Suppose that every
$\psi_l(\tau)$ $(l=1,\ldots,k)$ is a conti\-nu\-ous 
nonrandom function on the interval $[t, T]$
and $\{\phi_j(x)\}_{j=0}^{\infty}$ is a complete orthonormal system  
of continuous functions in the space $L_2([t,T]).$ Then

\vspace{-2mm}
$$
J[\psi^{(k)}]_{T,t}=
\hbox{\vtop{\offinterlineskip\halign{
\hfil#\hfil\cr
{\rm l.i.m.}\cr
$\stackrel{}{{}_{p_1,\ldots,p_k\to \infty}}$\cr
}} }J[\psi^{(k)}]_{T,t}^{p_1,\ldots,p_k},
$$

\vspace{3mm}
$$
{\sf M}\left\{\left(J[\psi^{(k)}]_{T,t}-J[\psi^{(k)}]_{T,t}^{p_1,\ldots,p_k}
\right)^2\right\}\le 
$$

\begin{equation}
\label{xx}
\le k!
\left(\int\limits_{[t,T]^k}
K^2(t_1,\ldots,t_k)dt_1\ldots dt_k-\sum_{j_1=0}^{p_1}\ldots
\sum_{j_k=0}^{p_k}C^2_{j_k\ldots j_1}\right),
\end{equation}

\newpage
\noindent
where
$$
J[\psi^{(k)}]_{T,t}^{p_1,\ldots,p_k}=
\sum_{j_1=0}^{p_1}\ldots\sum_{j_k=0}^{p_k}
C_{j_k\ldots j_1}\Biggl(
\prod_{l=1}^k\zeta_{j_l}^{(i_l)}-
\Biggr.
$$

\begin{equation}
\label{tyyy}
-\Biggl.
\hbox{\vtop{\offinterlineskip\halign{
\hfil#\hfil\cr
{\rm l.i.m.}\cr
$\stackrel{}{{}_{N\to \infty}}$\cr
}} }\sum_{(l_1,\ldots,l_k)\in {\rm G}_k}
\phi_{j_{1}}(\tau_{l_1})
\Delta{\bf w}_{\tau_{l_1}}^{(i_1)}\ldots
\phi_{j_{k}}(\tau_{l_k})
\Delta{\bf w}_{\tau_{l_k}}^{(i_k)}\Biggr)
\end{equation}

\vspace{2mm}
\noindent
and
$$
{\rm G}_k={\rm H}_k\backslash{\rm L}_k,\ \ \
{\rm H}_k=\{(l_1,\ldots,l_k):\ l_1,\ldots,l_k=0,\ 1,\ldots,N-1\},
$$
$$
{\rm L}_k=\{(l_1,\ldots,l_k):\ l_1,\ldots,l_k=0,\ 1,\ldots,N-1;\
l_g\ne l_r\ (g\ne r);\ g, r=1,\ldots,k\},
$$

\vspace{2mm}
\noindent
${\rm l.i.m.}$ is a limit in the mean-square sense{\rm ,} 
$i_1,\ldots,i_k=0,1,\ldots,m,$ 
\begin{equation}
\label{rr23}
\zeta_{j}^{(i)}=
\int\limits_t^T \phi_{j}(s) d{\bf w}_s^{(i)}
\end{equation} 
are independent standard Gaussian random variables
for various
$i$ or $j$ {\rm(}$ if~ i\ne 0${\rm),}
$C_{j_k\ldots j_1}$ is the Fourier coefficient {\rm(\ref{ppppa}),}
$\Delta{\bf w}_{\tau_{j}}^{(i)}=
{\bf w}_{\tau_{j+1}}^{(i)}-{\bf w}_{\tau_{j}}^{(i)}$
$(i=0,\ 1,\ldots,m),$\
$\left\{\tau_{j}\right\}_{j=0}^{N}$ 
is the discretization {\rm (\ref{1111})}, the estimate {\rm (\ref{xx})}
is valid for
$T-t\in (0, \infty)$ and $i_1,\ldots,i_k=1,\ldots,m$ or
$T-t\in (0, 1)$ and $i_1,\ldots,i_k=0, 1,\ldots,m$.
}

Note that in 
\cite{Kuznetsov3}-\cite{Kuznetsov5}, 
\cite{Kuznetsov6}, \cite{Kuznetsov7}, \cite{Kuznetsov8},
\cite{Kuznetsov16}
the version of Theorem 1 for
systems of Haar and Rademacher--Walsh 
functions has been considered. 
Some modifications of Theorem 1 for another types of iterated
stochastic integrals 
as well as for complete orthonormal with
weight $r(t_1)\ldots r(t_k)\ge 0$ systems of functions 
in the space $L_2([t, T]^k)$
can be found in 
\cite{Kuznetsov3}-\cite{Kuznetsov5},
\cite{Kuznetsov6},
\cite{Kuznetsov7}, \cite{Kuznetsov8}, 
\cite{Kuznetsov16}.
Application of Theorem 1 and Theorem 4 (see below) to
the mean-square approximation of iterated stochastic integrals
with respect to the infinite-dimensional $Q$-Wiener process
is presented in \cite{Kuznetsov13a}, \cite{Kuznetsov16} (Chapter 7),
\cite{arxiv-20}, \cite{arxiv-21}.

Obtain transformed particular cases of Theorem 1 for 
$k=1,\ldots,5$ \cite{Kuznetsov3}-\cite{Kuznetsov5},
\cite{Kuznetsov6},
\cite{Kuznetsov7}-\cite{Kuznetsov13a}, \cite{Kuznetsov15}, 
\cite{Kuznetsov16}

\vspace{-4mm}
\begin{equation}
\label{a1}
J[\psi^{(1)}]_{T,t}
=\hbox{\vtop{\offinterlineskip\halign{
\hfil#\hfil\cr
{\rm l.i.m.}\cr
$\stackrel{}{{}_{p_1\to \infty}}$\cr
}} }\sum_{j_1=0}^{p_1}
C_{j_1}\zeta_{j_1}^{(i_1)},
\end{equation}

\begin{equation}
\label{a2}
J[\psi^{(2)}]_{T,t}
=\hbox{\vtop{\offinterlineskip\halign{
\hfil#\hfil\cr
{\rm l.i.m.}\cr
$\stackrel{}{{}_{p_1,p_2\to \infty}}$\cr
}} }\sum_{j_1=0}^{p_1}\sum_{j_2=0}^{p_2}
C_{j_2j_1}\Biggl(\zeta_{j_1}^{(i_1)}\zeta_{j_2}^{(i_2)}
-{\bf 1}_{\{i_1=i_2\ne 0\}}
{\bf 1}_{\{j_1=j_2\}}\Biggr),
\end{equation}

$$
J[\psi^{(3)}]_{T,t}=
\hbox{\vtop{\offinterlineskip\halign{
\hfil#\hfil\cr
{\rm l.i.m.}\cr
$\stackrel{}{{}_{p_1,p_2,p_3\to \infty}}$\cr
}} }\sum_{j_1=0}^{p_1}\sum_{j_2=0}^{p_2}\sum_{j_3=0}^{p_3}
C_{j_3j_2j_1}\Biggl(
\zeta_{j_1}^{(i_1)}\zeta_{j_2}^{(i_2)}\zeta_{j_3}^{(i_3)}
-\Biggr.
$$
\begin{equation}
\label{a3}
-{\bf 1}_{\{i_1=i_2\ne 0\}}
{\bf 1}_{\{j_1=j_2\}}
\zeta_{j_3}^{(i_3)}
-{\bf 1}_{\{i_2=i_3\ne 0\}}
{\bf 1}_{\{j_2=j_3\}}
\zeta_{j_1}^{(i_1)}
\Biggl.-{\bf 1}_{\{i_1=i_3\ne 0\}}
{\bf 1}_{\{j_1=j_3\}}
\zeta_{j_2}^{(i_2)}\Biggr),
\end{equation}

\vspace{4mm}
$$
J[\psi^{(4)}]_{T,t}
=
\hbox{\vtop{\offinterlineskip\halign{
\hfil#\hfil\cr
{\rm l.i.m.}\cr
$\stackrel{}{{}_{p_1,\ldots,p_4\to \infty}}$\cr
}} }\sum_{j_1=0}^{p_1}\ldots\sum_{j_4=0}^{p_4}
C_{j_4\ldots j_1}\Biggl(
\prod_{l=1}^4\zeta_{j_l}^{(i_l)}
\Biggr.
-
$$
$$
-
{\bf 1}_{\{i_1=i_2\ne 0\}}
{\bf 1}_{\{j_1=j_2\}}
\zeta_{j_3}^{(i_3)}
\zeta_{j_4}^{(i_4)}
-
{\bf 1}_{\{i_1=i_3\ne 0\}}
{\bf 1}_{\{j_1=j_3\}}
\zeta_{j_2}^{(i_2)}
\zeta_{j_4}^{(i_4)}-
$$
$$
-
{\bf 1}_{\{i_1=i_4\ne 0\}}
{\bf 1}_{\{j_1=j_4\}}
\zeta_{j_2}^{(i_2)}
\zeta_{j_3}^{(i_3)}
-
{\bf 1}_{\{i_2=i_3\ne 0\}}
{\bf 1}_{\{j_2=j_3\}}
\zeta_{j_1}^{(i_1)}
\zeta_{j_4}^{(i_4)}-
$$
$$
-
{\bf 1}_{\{i_2=i_4\ne 0\}}
{\bf 1}_{\{j_2=j_4\}}
\zeta_{j_1}^{(i_1)}
\zeta_{j_3}^{(i_3)}
-
{\bf 1}_{\{i_3=i_4\ne 0\}}
{\bf 1}_{\{j_3=j_4\}}
\zeta_{j_1}^{(i_1)}
\zeta_{j_2}^{(i_2)}+
$$
$$
+
{\bf 1}_{\{i_1=i_2\ne 0\}}
{\bf 1}_{\{j_1=j_2\}}
{\bf 1}_{\{i_3=i_4\ne 0\}}
{\bf 1}_{\{j_3=j_4\}}
+
$$
$$
+
{\bf 1}_{\{i_1=i_3\ne 0\}}
{\bf 1}_{\{j_1=j_3\}}
{\bf 1}_{\{i_2=i_4\ne 0\}}
{\bf 1}_{\{j_2=j_4\}}+
$$
\begin{equation}
\label{a4}
+\Biggl.
{\bf 1}_{\{i_1=i_4\ne 0\}}
{\bf 1}_{\{j_1=j_4\}}
{\bf 1}_{\{i_2=i_3\ne 0\}}
{\bf 1}_{\{j_2=j_3\}}\Biggr),
\end{equation}

\vspace{5mm}

$$
J[\psi^{(5)}]_{T,t}
=\hbox{\vtop{\offinterlineskip\halign{
\hfil#\hfil\cr
{\rm l.i.m.}\cr
$\stackrel{}{{}_{p_1,\ldots,p_5\to \infty}}$\cr
}} }\sum_{j_1=0}^{p_1}\ldots\sum_{j_5=0}^{p_5}
C_{j_5\ldots j_1}\Biggl(
\prod_{l=1}^5\zeta_{j_l}^{(i_l)}
-\Biggr.
$$
$$
-
{\bf 1}_{\{i_1=i_2\ne 0\}}
{\bf 1}_{\{j_1=j_2\}}
\zeta_{j_3}^{(i_3)}
\zeta_{j_4}^{(i_4)}
\zeta_{j_5}^{(i_5)}-
{\bf 1}_{\{i_1=i_3\ne 0\}}
{\bf 1}_{\{j_1=j_3\}}
\zeta_{j_2}^{(i_2)}
\zeta_{j_4}^{(i_4)}
\zeta_{j_5}^{(i_5)}-
$$
$$
-
{\bf 1}_{\{i_1=i_4\ne 0\}}
{\bf 1}_{\{j_1=j_4\}}
\zeta_{j_2}^{(i_2)}
\zeta_{j_3}^{(i_3)}
\zeta_{j_5}^{(i_5)}-
{\bf 1}_{\{i_1=i_5\ne 0\}}
{\bf 1}_{\{j_1=j_5\}}
\zeta_{j_2}^{(i_2)}
\zeta_{j_3}^{(i_3)}
\zeta_{j_4}^{(i_4)}-
$$
$$
-
{\bf 1}_{\{i_2=i_3\ne 0\}}
{\bf 1}_{\{j_2=j_3\}}
\zeta_{j_1}^{(i_1)}
\zeta_{j_4}^{(i_4)}
\zeta_{j_5}^{(i_5)}-
{\bf 1}_{\{i_2=i_4\ne 0\}}
{\bf 1}_{\{j_2=j_4\}}
\zeta_{j_1}^{(i_1)}
\zeta_{j_3}^{(i_3)}
\zeta_{j_5}^{(i_5)}-
$$
$$
-
{\bf 1}_{\{i_2=i_5\ne 0\}}
{\bf 1}_{\{j_2=j_5\}}
\zeta_{j_1}^{(i_1)}
\zeta_{j_3}^{(i_3)}
\zeta_{j_4}^{(i_4)}
-{\bf 1}_{\{i_3=i_4\ne 0\}}
{\bf 1}_{\{j_3=j_4\}}
\zeta_{j_1}^{(i_1)}
\zeta_{j_2}^{(i_2)}
\zeta_{j_5}^{(i_5)}-
$$
$$
-
{\bf 1}_{\{i_3=i_5\ne 0\}}
{\bf 1}_{\{j_3=j_5\}}
\zeta_{j_1}^{(i_1)}
\zeta_{j_2}^{(i_2)}
\zeta_{j_4}^{(i_4)}
-{\bf 1}_{\{i_4=i_5\ne 0\}}
{\bf 1}_{\{j_4=j_5\}}
\zeta_{j_1}^{(i_1)}
\zeta_{j_2}^{(i_2)}
\zeta_{j_3}^{(i_3)}+
$$
$$
+
{\bf 1}_{\{i_1=i_2\ne 0\}}
{\bf 1}_{\{j_1=j_2\}}
{\bf 1}_{\{i_3=i_4\ne 0\}}
{\bf 1}_{\{j_3=j_4\}}\zeta_{j_5}^{(i_5)}+
{\bf 1}_{\{i_1=i_2\ne 0\}}
{\bf 1}_{\{j_1=j_2\}}
{\bf 1}_{\{i_3=i_5\ne 0\}}
{\bf 1}_{\{j_3=j_5\}}\zeta_{j_4}^{(i_4)}+
$$
$$
+
{\bf 1}_{\{i_1=i_2\ne 0\}}
{\bf 1}_{\{j_1=j_2\}}
{\bf 1}_{\{i_4=i_5\ne 0\}}
{\bf 1}_{\{j_4=j_5\}}\zeta_{j_3}^{(i_3)}+
{\bf 1}_{\{i_1=i_3\ne 0\}}
{\bf 1}_{\{j_1=j_3\}}
{\bf 1}_{\{i_2=i_4\ne 0\}}
{\bf 1}_{\{j_2=j_4\}}\zeta_{j_5}^{(i_5)}+
$$
$$
+
{\bf 1}_{\{i_1=i_3\ne 0\}}
{\bf 1}_{\{j_1=j_3\}}
{\bf 1}_{\{i_2=i_5\ne 0\}}
{\bf 1}_{\{j_2=j_5\}}\zeta_{j_4}^{(i_4)}+
{\bf 1}_{\{i_1=i_3\ne 0\}}
{\bf 1}_{\{j_1=j_3\}}
{\bf 1}_{\{i_4=i_5\ne 0\}}
{\bf 1}_{\{j_4=j_5\}}\zeta_{j_2}^{(i_2)}+
$$
$$
+
{\bf 1}_{\{i_1=i_4\ne 0\}}
{\bf 1}_{\{j_1=j_4\}}
{\bf 1}_{\{i_2=i_3\ne 0\}}
{\bf 1}_{\{j_2=j_3\}}\zeta_{j_5}^{(i_5)}+
{\bf 1}_{\{i_1=i_4\ne 0\}}
{\bf 1}_{\{j_1=j_4\}}
{\bf 1}_{\{i_2=i_5\ne 0\}}
{\bf 1}_{\{j_2=j_5\}}\zeta_{j_3}^{(i_3)}+
$$
$$
+
{\bf 1}_{\{i_1=i_4\ne 0\}}
{\bf 1}_{\{j_1=j_4\}}
{\bf 1}_{\{i_3=i_5\ne 0\}}
{\bf 1}_{\{j_3=j_5\}}\zeta_{j_2}^{(i_2)}+
{\bf 1}_{\{i_1=i_5\ne 0\}}
{\bf 1}_{\{j_1=j_5\}}
{\bf 1}_{\{i_2=i_3\ne 0\}}
{\bf 1}_{\{j_2=j_3\}}\zeta_{j_4}^{(i_4)}+
$$
$$
+
{\bf 1}_{\{i_1=i_5\ne 0\}}
{\bf 1}_{\{j_1=j_5\}}
{\bf 1}_{\{i_2=i_4\ne 0\}}
{\bf 1}_{\{j_2=j_4\}}\zeta_{j_3}^{(i_3)}+
{\bf 1}_{\{i_1=i_5\ne 0\}}
{\bf 1}_{\{j_1=j_5\}}
{\bf 1}_{\{i_3=i_4\ne 0\}}
{\bf 1}_{\{j_3=j_4\}}\zeta_{j_2}^{(i_2)}+
$$
$$
+
{\bf 1}_{\{i_2=i_3\ne 0\}}
{\bf 1}_{\{j_2=j_3\}}
{\bf 1}_{\{i_4=i_5\ne 0\}}
{\bf 1}_{\{j_4=j_5\}}\zeta_{j_1}^{(i_1)}+
{\bf 1}_{\{i_2=i_4\ne 0\}}
{\bf 1}_{\{j_2=j_4\}}
{\bf 1}_{\{i_3=i_5\ne 0\}}
{\bf 1}_{\{j_3=j_5\}}\zeta_{j_1}^{(i_1)}+
$$
\begin{equation}
\label{a5}
+\Biggl.
{\bf 1}_{\{i_2=i_5\ne 0\}}
{\bf 1}_{\{j_2=j_5\}}
{\bf 1}_{\{i_3=i_4\ne 0\}}
{\bf 1}_{\{j_3=j_4\}}\zeta_{j_1}^{(i_1)}\Biggr),
\end{equation}

\noindent
where ${\bf 1}_A$ is the indicator of the set $A$.

Let us consider the generalization of the formulas (\ref{a1})--(\ref{a5}) 
for the case of an arbitrary $k$ ($k\in\mathbb{N}$). 

{\bf Theorem 2} \cite{Kuznetsov4a} (2009), \cite{Kuznetsov5},
\cite{Kuznetsov6}, \cite{Kuznetsov7}, \cite{Kuznetsov8},
\cite{Kuznetsov13a},
\cite{Kuznetsov16}.
{\it Under the conditions of Theorem {\rm 1} 
the following expansion 
$$
J[\psi^{(k)}]_{T,t}^{(i_1\ldots i_k)}=
\hbox{\vtop{\offinterlineskip\halign{
\hfil#\hfil\cr
{\rm l.i.m.}\cr
$\stackrel{}{{}_{p_1,\ldots,p_k\to \infty}}$\cr
}} }
\sum\limits_{j_1=0}^{p_1}\ldots
\sum\limits_{j_k=0}^{p_k}
C_{j_k\ldots j_1}\Biggl(
\prod_{l=1}^k\zeta_{j_l}^{(i_l)}+\sum\limits_{r=1}^{[k/2]}
(-1)^r \times
\Biggr.
$$

\vspace{-4mm}
\begin{equation}
\label{leto6000}
\times
\sum_{\stackrel{(\{\{g_1, g_2\}, \ldots, 
\{g_{2r-1}, g_{2r}\}\}, \{q_1, \ldots, q_{k-2r}\})}
{{}_{\{g_1, g_2, \ldots, 
g_{2r-1}, g_{2r}, q_1, \ldots, q_{k-2r}\}=\{1, 2, \ldots, k\}}}}
\prod\limits_{s=1}^r
{\bf 1}_{\{i_{g_{{}_{2s-1}}}=~i_{g_{{}_{2s}}}\ne 0\}}
\Biggl.{\bf 1}_{\{j_{g_{{}_{2s-1}}}=~j_{g_{{}_{2s}}}\}}
\prod_{l=1}^{k-2r}\zeta_{j_{q_l}}^{(i_{q_l})}\Biggr)
\end{equation}

\vspace{2mm}
\noindent
converging in the mean-square sense is valid, 
where $[\cdot]$ is an integer part of a real number,
$$
\sum_{\stackrel{(\{\{g_1, g_2\}, \ldots, 
\{g_{2r-1}, g_{2r}\}\}, \{q_1, \ldots, q_{k-2r}\})}
{{}_{\{g_1, g_2, \ldots, 
g_{2r-1}, g_{2r}, q_1, \ldots, q_{k-2r}\}=\{1, 2, \ldots, k\}}}}
$$

\noindent
means the sum with respect to all possible permutations of the set
$$
(\{\{g_1, g_2\}, \ldots, 
\{g_{2r-1}, g_{2r}\}\}, \{q_1, \ldots, q_{k-2r}\}),
$$

\noindent
where $\{g_1, g_2, \ldots, 
g_{2r-1}, g_{2r}, q_1, \ldots, q_{k-2r}\}=\{1, 2, \ldots, k\},$
braces mean an unordered 
set, and parentheses mean an ordered set{\rm ;} another notations
are the same as in Theorem {\rm 1}.
}

\vspace{2mm}

For further consideration, we need the following statement.
                    
{\bf Theorem 3}\ \cite{Kuznetsov16} (Sect.~1.1.9, 1.11, 1.12),\ \cite{arxiv-1} (Sect.~6, 15, 16).\
{\it Under the conditions of Theorem {\rm 1} 
the following estimate 

\newpage
\noindent
$$
{\sf M}\left\{\left(J[\psi^{(k)}]_{T,t}-J[\psi^{(k)}]_{T,t}^{p_1,\ldots,p_k}
\right)^{2n}\right\}\le 
$$

$$
\le
(k!)^{n} (2n-1)^{nk}\ \times
$$

\vspace{-2mm}
\begin{equation}
\label{dima2ye100}
\times\ 
\left(
\int\limits_{[t,T]^k}
K^2(t_1,\ldots,t_k)
dt_1\ldots dt_k -\sum_{j_1=0}^{p_1}\ldots
\sum_{j_k=0}^{p_k}C^2_{j_k\ldots j_1}
\right)^n
\end{equation}

\vspace{3mm}
\noindent
is valid, where $n\in\mathbb{N};$ another notations
are the same as in Theorem {\rm 1}.}

Since according
to the Parseval's equality
$$
\int\limits_{[t,T]^k}
K^2(t_1,\ldots,t_k)
dt_1\ldots dt_k -\sum_{j_1=0}^{p_1}\ldots
\sum_{j_k=0}^{p_k}C^2_{j_k\ldots j_1}\  \to \  0
$$

\noindent
if $p_1,\ldots,p_k\to\infty$,
then the inequality (\ref{dima2ye100})
means that the expansions of 
iterated It\^{o} stochastic integrals in Theorem 1
converge in the 
mean
of degree $2n$ ($n\in \mathbb{N}$).

Let us consider the generalization of Theorems 1--3 for the case
of an arbitrary complete orthonormal systems  
of functions in the space $L_2([t,T])$ 
and $\psi_1(\tau),\ldots,\psi_k(\tau)\in L_2([t, T]).$

{\bf Theorem~4}\ \cite{Kuznetsov16} (Sect.~1.1.9, 1.11, 1.12),\ \cite{arxiv-1} (Sect.~6, 15, 16).\
{\it Suppose that
$\psi_1(\tau),$ $\ldots,$ $\psi_k(\tau)$ $\in L_2([t, T])$ and
$\{\phi_j(x)\}_{j=0}^{\infty}$ is an arbitrary complete orthonormal system  
of functions in the space $L_2([t,T]).$
Then 

\vspace{-1mm}
\begin{equation}
\label{yyy1}
J[\psi^{(k)}]_{T,t}=
\hbox{\vtop{\offinterlineskip\halign{
\hfil#\hfil\cr
{\rm l.i.m.}\cr
$\stackrel{}{{}_{p_1,\ldots,p_k\to \infty}}$\cr
}} }J[\psi^{(k)}]_{T,t}^{p_1,\ldots,p_k},
\end{equation}

\vspace{3mm}
$$
{\sf M}\left\{\left(J[\psi^{(k)}]_{T,t}-J[\psi^{(k)}]_{T,t}^{p_1,\ldots,p_k}
\right)^2\right\}\le 
$$

$$
\le k!
\left(\int\limits_{[t,T]^k}
K^2(t_1,\ldots,t_k)dt_1\ldots dt_k-\sum_{j_1=0}^{p_1}\ldots
\sum_{j_k=0}^{p_k}C^2_{j_k\ldots j_1}\right),
$$

$$
{\sf M}\left\{\left(J[\psi^{(k)}]_{T,t}-J[\psi^{(k)}]_{T,t}^{p_1,\ldots,p_k}
\right)^{2n}\right\}\le 
$$

$$
\le
(k!)^{n} (2n-1)^{nk}\ \times
$$

\vspace{-2mm}
$$
\times\ 
\left(
\int\limits_{[t,T]^k}
K^2(t_1,\ldots,t_k)
dt_1\ldots dt_k -\sum_{j_1=0}^{p_1}\ldots
\sum_{j_k=0}^{p_k}C^2_{j_k\ldots j_1}
\right)^n,
$$

\vspace{4mm}
\noindent
where $n\in\mathbb{N},$
$$
J[\psi^{(k)}]_{T,t}^{p_1,\ldots,p_k}=
\sum\limits_{j_1=0}^{p_1}\ldots
\sum\limits_{j_k=0}^{p_k}
C_{j_k\ldots j_1}\Biggl(
\prod_{l=1}^k\zeta_{j_l}^{(i_l)}+\sum\limits_{r=1}^{[k/2]}
(-1)^r \times
\Biggr.
$$

\vspace{-4mm}
\begin{equation}
\label{xxx333}
\times
\sum_{\stackrel{(\{\{g_1, g_2\}, \ldots, 
\{g_{2r-1}, g_{2r}\}\}, \{q_1, \ldots, q_{k-2r}\})}
{{}_{\{g_1, g_2, \ldots, 
g_{2r-1}, g_{2r}, q_1, \ldots, q_{k-2r}\}=\{1, 2, \ldots, k\}}}}
\prod\limits_{s=1}^r
{\bf 1}_{\{i_{g_{{}_{2s-1}}}=~i_{g_{{}_{2s}}}\ne 0\}}
\Biggl.{\bf 1}_{\{j_{g_{{}_{2s-1}}}=~j_{g_{{}_{2s}}}\}}
\prod_{l=1}^{k-2r}\zeta_{j_{q_l}}^{(i_{q_l})}\Biggr),
\end{equation}

\vspace{2mm}
\noindent
where $[x]$ is an integer part of a real number $x;$
another notations are the same as in Theorems~{\rm 1--3}.}

It should be noted that an analogue of the expansion (\ref{yyy1})
under the conditions of Theorem 4
was considered 
in \cite{Rybakov1000}. 
Note that we use another notations 
\cite{Kuznetsov16} (Sect.~1.11), \cite{arxiv-1} (Sect.~15)
in comparison with \cite{Rybakov1000}.
Moreover, the proof of an analogue of (\ref{yyy1})
from \cite{Rybakov1000} is somewhat different from the proof given in 
\cite{Kuznetsov16} (Sect.~1.11), \cite{arxiv-1} (Sect.~15).

Also note the following theorem.

{\bf Theorem 5} \cite{Kuznetsov16} (Sect.~1.12), \cite{arxiv-3} (Sect.~6).
{\it Suppose that $\{\phi_j(x)\}_{j=0}^{\infty}$ 
is an arbitrary complete orthonormal system  
of functions in the space $L_2([t,T])$ and
$\psi_1(\tau),\ldots,\psi_k(\tau)\in L_2([t, T]),$  $i_1,\ldots, i_k=1,\ldots,m$.
Then

\vspace{-2mm}
$$
{\sf M}\left\{\left(J[\psi^{(k)}]_{T,t}-J[\psi^{(k)}]_{T,t}^{p,\ldots,p}
\right)^{2}\right\}=
\int\limits_{[t,T]^k}
K^2(t_1,\ldots,t_k)dt_1\ldots dt_k- 
$$
$$
-\sum_{j_1,\ldots, j_k=0}^{p}
C_{j_k\ldots j_1}
{\sf M}\left\{J[\psi^{(k)}]_{T,t}
\sum\limits_{(j_1,\ldots,j_k)}
\int\limits_t^T \phi_{j_k}(t_k)
\ldots
\int\limits_t^{t_{2}}\phi_{j_{1}}(t_{1})
d{\bf w}_{t_1}^{(i_1)}\ldots
d{\bf w}_{t_k}^{(i_k)}\right\},
$$

\vspace{3mm}
\noindent
where $i_1,\ldots,i_k = 1,\ldots,m;$ the value $J[\psi^{(k)}]_{T,t}^{p,\ldots,p}$
is defined by {\rm (\ref{xxx333})} {\rm ($p_1=\ldots=p_k=p$);}
the expression 
$$
\sum\limits_{(j_1,\ldots,j_k)}
$$

\noindent
means the sum with respect to all
possible permutations 
$(j_1,\ldots,j_k)$. At the same time if 
$j_r$ swapped with $j_q$ in the permutation $(j_1,\ldots,j_k),$
then $i_r$ swapped with $i_q$ in the permutation
$(i_1,\ldots,i_k);$ 
another notations are the same as in Theorems {\rm 1, 2.}
}

Let us consider the following iterated It\^{o} stochastic 
integrals from the Taylor--It\^{o} expansion \cite{Kloeden2}
\begin{equation}
\label{k1000}
J_{(\lambda_1\ldots \lambda_k)T,t}^{(i_1\ldots i_k)}
=\int\limits_t^T\ldots \int\limits_t^{t_{2}}
d{\bf w}_{t_1}^{(i_1)}\ldots
d{\bf w}_{t_k}^{(i_k)},
\end{equation}

\noindent
where $i_1,\ldots, i_k=0, 1,\dots,m,$
$\lambda_l=1$ if $i_l=1,\ldots,m$ and 
$\lambda_l=0$ if $i_l=0$ $(l=1,\ldots,k)$.
Remind that 
${\bf w}_{\tau}^{(i)},$ $i=1,\ldots,m$ are independent
standard Wiener processes and
${\bf w}_{\tau}^{(0)}=\tau.$

For example,
using Theorems 1, 4 (see (\ref{a1})-(\ref{a3})) and 
complete orthonormal system of Legendre polynomials in 
the space $L_2([t,T])$ 
we obtain the following approximations of 
the iterated 
It\^{o} stochastic integrals (\ref{k1000})
\cite{Kuznetsov3}-\cite{Kuznetsov5},
\cite{Kuznetsov6},
\cite{Kuznetsov7}-\cite{Kuznetsov13a}, \cite{Kuznetsov15},
\cite{Kuznetsov16}
(also see
early publications \cite{Kuznetsov1}, \cite{Kuznetsov2})

\vspace{-1mm}
\begin{equation}
\label{opp0}
J_{(1)T,t}^{(i_1)}=\sqrt{T-t}\zeta_0^{(i_1)},
\end{equation}

\begin{equation}
\label{opp1}
J_{(01)T,t}^{(0 i_1)}=\frac{(T-t)^{3/2}}{2}\biggl(\zeta_0^{(i_1)}+
\frac{1}{\sqrt{3}}\zeta_1^{(i_1)}\biggr),
\end{equation}

\begin{equation}
\label{opp2}
J_{(10)T,t}^{(i_1 0)}=\frac{(T-t)^{3/2}}{2}\biggl(\zeta_0^{(i_1)}-
\frac{1}{\sqrt{3}}\zeta_1^{(i_1)}\biggr),
\end{equation}

\vspace{-1mm}
\begin{equation}
\label{kr00}
J_{(11)T,t}^{(i_1 i_2)q}=
\frac{T-t}{2}\Biggl(\zeta_0^{(i_1)}\zeta_0^{(i_2)}+
\sum_{i=1}^{q}
\frac{1}{\sqrt{4i^2-1}}\biggl(
\zeta_{i-1}^{(i_1)}\zeta_{i}^{(i_2)}-
\zeta_i^{(i_1)}\zeta_{i-1}^{(i_2)}\biggr)-{\bf 1}_{\{i_1=i_2\}}
\Biggr),
\end{equation}

\vspace{2mm}
$$
J_{(11)T,t}^{(i_1 i_1)}
=\frac{1}{2}(T-t)\biggl(
\left(\zeta_0^{(i_1)}\right)^2-1\biggr),
$$

$$
J_{(111)T,t}^{(i_1 i_2 i_3)p}=
\sum_{j_1,j_2,j_3=0}^{p}
C_{j_3j_2j_1}
\Biggl(
\zeta_{j_1}^{(i_1)}\zeta_{j_2}^{(i_2)}\zeta_{j_3}^{(i_3)}
\Biggr.-{\bf 1}_{\{i_1=i_2\}}
{\bf 1}_{\{j_1=j_2\}}
\zeta_{j_3}^{(i_3)}-
$$
\begin{equation}
\label{kr1}
\Biggl.-{\bf 1}_{\{i_2=i_3\}}
{\bf 1}_{\{j_2=j_3\}}
\zeta_{j_1}^{(i_1)}-
{\bf 1}_{\{i_1=i_3\}}
{\bf 1}_{\{j_1=j_3\}}
\zeta_{j_2}^{(i_2)}\Biggr),
\end{equation}

\vspace{4mm}
$$
J_{(111)T,t}^{(i_1 i_1 i_1)}
=\frac{1}{6}(T-t)^{3/2}\biggl(
\left(\zeta_0^{(i_1)}\right)^3-3
\zeta_0^{(i_1)}\biggr),
$$

\vspace{4mm}
\noindent
where

\vspace{-3mm}
$$
C_{j_3j_2j_1}=\frac{\sqrt{(2j_1+1)(2j_2+1)(2j_3+1)}(T-t)^{3/2}}{8}\bar
C_{j_3j_2j_1},
$$

$$
\bar C_{j_3j_2j_1}=\int\limits_{-1}^{1}P_{j_3}(z)
\int\limits_{-1}^{z}P_{j_2}(y)
\int\limits_{-1}^{y}
P_{j_1}(x)dx dy dz,
$$

\vspace{3mm}
\noindent
where the 
Gaussian random variable $\zeta_{j}^{(i)}$ (if $i\ne 0$) is defined by 
(\ref{rr23}) and 
$P_{j}(x)$ $(j=0, 1, 2,\ldots )$ is the 
Legendre polynomial \cite{Suetin}.

Note that formula (\ref{kr00}) has been 
obtained for the first time in \cite{Kuznetsov1} (1997).
For pairwise different $i_1,i_2,i_3=1,\ldots,m$
we have \cite{Kuznetsov1}, \cite{Kuznetsov2},
\cite{Kuznetsov3}-\cite{Kuznetsov5},
\cite{Kuznetsov6},
\cite{Kuznetsov7}-\cite{Kuznetsov13a}, \cite{Kuznetsov15}

\vspace{-3mm}
\begin{equation}
\label{909}
{\sf M}\left\{\left(J_{(11)T,t}^{(i_1 i_2)}-J_{(11)T,t}^{(i_1 i_2)q}
\right)^2\right\}
=\frac{(T-t)^2}{2}\left(\frac{1}{2}-\sum_{i=1}^{q}
\frac{1}{4i^2-1}\right),
\end{equation}

\begin{equation}
\label{zzz}
{\sf M}\left\{\left(J_{(111)T,t}^{(i_1 i_2 i_3)}-J_{(111)T,t}^{(i_1 i_2 i_3)p}
\right)^2\right\}=
\frac{(T-t)^{3}}{6}-\sum_{j_1,j_2,j_3=0}^{p}
C_{j_3j_2j_1}^2.
\end{equation}

\vspace{3mm}

The problem of the exact calculation of the mean-square 
error of approximation in Theorems 1, 4 is solved completely 
for an arbitrary $k$ $(k\in \mathbb{N})$ and
any possible combinations of the numbers $i_1,\ldots,i_k=1,\ldots,m$
in Theorem 5
(also see \cite{Kuznetsov8}, \cite{Kuznetsov16}, \cite{arxiv-3}).

\section{Convergence With Probability 1
of Expansions of Iterated
It\^{o} Stochastic Integrals of Multiplicity $k$ $(k\in\mathbb{N})$
in Theorems 1, 2}

Let us address now to the convergence (w.~p.~1)
in Theorem 1. As we mentioned above this question has been studied
for simplest iterated It\^{o} stochastic integrals
of multiplicities 1 and 2 in 
\cite{Kuznetsov4}-\cite{Kuznetsov5},
\cite{Kuznetsov6},
\cite{Kuznetsov7}, \cite{Kuznetsov8}, \cite{Kuznetsov16}.

In this section, we formulate and prove the general result on 
convergence w. p. 1 of expansions 
of iterated It\^{o} stochastic integrals in Theorems 1, 2
for the case of multiplicity $k$ $(k\in\mathbb{N})$
for these integrals.

{\bf Theorem 6.} {\it Let 
$\psi_l(\tau)$ $(l=1,\ldots, k)$ are 
continuously differentiable nonrandom functions on the interval
$[t, T]$ and $\{\phi_j(x)\}_{j=0}^{\infty}$ is a complete
orthonormal system of Legendre polynomials or 
trigonometric functions in the space $L_2([t, T]).$
Then $J[\psi^{(k)}]_{T,t}^{p,\ldots,p}\ \to \ J[\psi^{(k)}]_{T,t}$
if $p\to \infty$ w.~p.~{\rm 1},
where $J[\psi^{(k)}]_{T,t}^{p,\ldots,p}$
is defined as the right-hand side of {\rm (\ref{leto6000})}
before passing to the limit
for the case $p_1=\ldots=p_k=p,$ i.e. {\rm (}see Theorem {\rm 2)}
$$
J[\psi^{(k)}]_{T,t}^{p,\ldots,p}=
\sum\limits_{j_1,\ldots, j_k=0}^{p}
C_{j_k\ldots j_1}\Biggl(
\prod_{l=1}^k\zeta_{j_l}^{(i_l)}+\sum\limits_{r=1}^{[k/2]}
(-1)^r \times
\Biggr.
$$

\vspace{-4mm}
$$
\times
\sum_{\stackrel{(\{\{g_1, g_2\}, \ldots, 
\{g_{2r-1}, g_{2r}\}\}, \{q_1, \ldots, q_{k-2r}\})}
{{}_{\{g_1, g_2, \ldots, 
g_{2r-1}, g_{2r}, q_1, \ldots, q_{k-2r}\}=\{1, 2, \ldots, k\}}}}
\prod\limits_{s=1}^r
{\bf 1}_{\{i_{g_{{}_{2s-1}}}=~i_{g_{{}_{2s}}}\ne 0\}}
\Biggl.{\bf 1}_{\{j_{g_{{}_{2s-1}}}=~j_{g_{{}_{2s}}}\}}
\prod_{l=1}^{k-2r}\zeta_{j_{q_l}}^{(i_{q_l})}\Biggr),
$$

\vspace{2mm}
\noindent
where $i_1,\ldots,i_k=1,\ldots,m,$ another notations are the same as in 
Theorems {\rm 1, 2.}}

{\bf Proof.} Let us consider the Parseval equality
\begin{equation}
\label{par1}
\int\limits_{[t,T]^k}K^2(t_1,\ldots,t_k)dt_1\ldots dt_k=
\lim\limits_{p_1,\ldots,p_k\to\infty}
\sum_{j_1=0}^{p_1}\ldots \sum_{j_k=0}^{p_k}
C_{j_k\ldots j_1}^2,
\end{equation}

\noindent
where
$$
K(t_1,\ldots,t_k)=
\begin{cases}
\psi_1(t_1)\ldots \psi_k(t_k),\ &t_1<\ldots<t_k\\
~\\
0,\ &\hbox{\rm otherwise}
\end{cases},
$$

\newpage
\noindent
where $t_1,\ldots,t_k\in [t, T]$ ($k\ge 2$) and 
$K(t_1)\equiv\psi_1(t_1)$ for $t_1\in[t, T],$ 
$$
C_{j_k\ldots j_1}=\int\limits_{[t,T]^k}
K(t_1,\ldots,t_k)\prod_{l=1}^{k}\phi_{j_l}(t_l)dt_1\ldots dt_k
$$
is the Fourier coefficient.

Taking into account the definitions of $K(t_1,\ldots,t_k)$
and $C_{j_k\ldots j_1}$, we obtain
\begin{equation}
\label{rrr1}
C_{j_k\ldots j_1}=
\int\limits_t^T
\phi_{j_k}(t_k)\psi_k(t_k)\ldots \int\limits_t^{t_2}
\phi_{j_1}(t_1)\psi_1(t_1)dt_1\ldots dt_k.
\end{equation}

Further, we denote
$$
\lim\limits_{p_1,\ldots,p_k\to\infty}
\sum_{j_1=0}^{p_1}\ldots \sum_{j_k=0}^{p_k}
C_{j_k\ldots j_1}^2\stackrel{\sf def}{=}
\sum_{j_1,\ldots,j_k=0}^{\infty}
C_{j_k\ldots j_1}^2.
$$

If $p_1=\ldots=p_k=p,$ then we also write
$$
\lim\limits_{p\to\infty}
\sum_{j_1=0}^{p}\ldots \sum_{j_k=0}^{p}
C_{j_k\ldots j_1}^2\stackrel{\sf def}{=}
\sum_{j_1,\ldots,j_k=0}^{\infty}
C_{j_k\ldots j_1}^2.
$$

From the other hand, for iterated limits we write
$$
\lim\limits_{p_1\to\infty}\ldots \lim\limits_{p_k\to\infty}
\sum_{j_1=0}^{p_1}\ldots \sum_{j_k=0}^{p_k}
C_{j_k\ldots j_1}^2\stackrel{\sf def}{=}
\sum_{j_1=0}^{\infty}\ldots
\sum_{j_k=0}^{\infty}
C_{j_k\ldots j_1}^2,
$$
$$
\lim\limits_{p_1\to\infty}\lim\limits_{p_2,\ldots,p_k\to\infty}
\sum_{j_1=0}^{p_1}\ldots \sum_{j_k=0}^{p_k}
C_{j_k\ldots j_1}^2\stackrel{\sf def}{=}
\sum_{j_1=0}^{\infty}
\sum_{j_2,\ldots,j_k=0}^{\infty}
C_{j_k\ldots j_1}^2
$$
and so on.

Using the Parseval equality and Lemma 2 (see Appendix) we obtain
$$
\int\limits_{[t,T]^k}K^2(t_1,\ldots,t_k)dt_1\ldots dt_k-
\sum_{j_1=0}^{p}\ldots \sum_{j_k=0}^{p}
C_{j_k\ldots j_1}^2=
$$
$$
=\sum_{j_1,\ldots,j_k=0}^{\infty}
C_{j_k\ldots j_1}^2-
\sum_{j_1=0}^{p}\ldots \sum_{j_k=0}^{p}
C_{j_k\ldots j_1}^2=
$$

\newpage
\noindent
$$
=
\sum_{j_1=0}^{\infty}\ldots \sum_{j_k=0}^{\infty}
C_{j_k\ldots j_1}^2-
\sum_{j_1=0}^{p}\ldots \sum_{j_k=0}^{p}
C_{j_k\ldots j_1}^2=
$$
$$
=\sum_{j_1=0}^{p}\sum_{j_2=0}^{\infty}\ldots \sum_{j_k=0}^{\infty}
C_{j_k\ldots j_1}^2+
\sum_{j_1=p+1}^{\infty}\sum_{j_2=0}^{\infty}\ldots \sum_{j_k=0}^{\infty}
C_{j_k\ldots j_1}^2-
\sum_{j_1=0}^{p}\ldots \sum_{j_k=0}^{p}
C_{j_k\ldots j_1}^2=
$$
$$
=\sum_{j_1=0}^{p}\sum_{j_2=0}^{p}\sum_{j_3=0}^{\infty}
\ldots \sum_{j_k=0}^{\infty}
C_{j_k\ldots j_1}^2+
\sum_{j_1=0}^{p}\sum_{j_2=p+1}^{\infty}
\sum_{j_3=0}^{\infty}
\ldots \sum_{j_k=0}^{\infty}+
$$
$$
+\sum_{j_1=p+1}^{\infty}\sum_{j_2=0}^{\infty}\ldots \sum_{j_k=0}^{\infty}
C_{j_k\ldots j_1}^2-
\sum_{j_1=0}^{p}\ldots \sum_{j_k=0}^{p}
C_{j_k\ldots j_1}^2=
$$
$$
=\ldots =
$$
$$
=\sum_{j_1=p+1}^{\infty}\sum_{j_2=0}^{\infty}\ldots \sum_{j_k=0}^{\infty}
C_{j_k\ldots j_1}^2+
\sum_{j_1=0}^p
\sum_{j_2=p+1}^{\infty}\sum_{j_2=0}^{\infty}\ldots \sum_{j_k=0}^{\infty}
C_{j_k\ldots j_1}^2+
$$
$$
+\sum_{j_1=0}^p\sum_{j_2=0}^p
\sum_{j_3=p+1}^{\infty}\sum_{j_4=0}^{\infty}\ldots \sum_{j_k=0}^{\infty}
C_{j_k\ldots j_1}^2+ \ldots +
\sum_{j_1=0}^p\ldots \sum_{j_{k-1}=0}^p
\sum_{j_k=p+1}^{\infty}C_{j_k\ldots j_1}^2\le
$$
$$
\le\sum_{j_1=p+1}^{\infty}\sum_{j_2=0}^{\infty}\ldots \sum_{j_k=0}^{\infty}
C_{j_k\ldots j_1}^2+
\sum_{j_1=0}^{\infty}
\sum_{j_2=p+1}^{\infty}\sum_{j_2=0}^{\infty}\ldots \sum_{j_k=0}^{\infty}
C_{j_k\ldots j_1}^2+
$$
$$
+\sum_{j_1=0}^{\infty}\sum_{j_2=0}^{\infty}
\sum_{j_3=p+1}^{\infty}\sum_{j_4=0}^{\infty}\ldots \sum_{j_k=0}^{\infty}
C_{j_k\ldots j_1}^2+ \ldots +
\sum_{j_1=0}^{\infty}\ldots \sum_{j_{k-1}=0}^{\infty}
\sum_{j_k=p+1}^{\infty}C_{j_k\ldots j_1}^2=
$$
\begin{equation}
\label{aaap}
=\sum\limits_{s=1}^k \left(\sum_{j_1=0}^{\infty}\ldots
\sum_{j_{s-1}=0}^{\infty}
\sum_{j_s=p+1}^{\infty}\sum_{j_{s+1}=0}^{\infty}\ldots \sum_{j_k=0}^{\infty}
C_{j_k\ldots j_1}^2\right).
\end{equation}

\vspace{2mm}

Note that deriving (\ref{aaap}) we use the following
$$
\sum_{j_1=0}^{p}\ldots
\sum_{j_{s-1}=0}^{p}
\sum_{j_s=p+1}^{\infty}\sum_{j_{s+1}=0}^{\infty}\ldots \sum_{j_k=0}^{\infty}
C_{j_k\ldots j_1}^2\le
$$
$$
\le
\sum_{j_1=0}^{m_1}\ldots
\sum_{j_{s-1}=0}^{m_{s-1}}
\sum_{j_s=p+1}^{\infty}\sum_{j_{s+1}=0}^{\infty}\ldots \sum_{j_k=0}^{\infty}
C_{j_k\ldots j_1}^2\le
$$
$$
\le
\lim\limits_{m_{s-1}\to\infty}
\sum_{j_1=0}^{m_1}\ldots
\sum_{j_{s-1}=0}^{m_{s-1}}
\sum_{j_s=p+1}^{\infty}\sum_{j_{s+1}=0}^{\infty}\ldots \sum_{j_k=0}^{\infty}
C_{j_k\ldots j_1}^2=
$$

\newpage
\noindent
$$
=
\sum_{j_1=0}^{m_1}\ldots
\sum_{j_{s-2}=0}^{m_{s-2}}\sum_{j_{s-1}=0}^{\infty}
\sum_{j_s=p+1}^{\infty}\sum_{j_{s+1}=0}^{\infty}\ldots \sum_{j_k=0}^{\infty}
C_{j_k\ldots j_1}^2\le
\ldots\le
$$
$$
\le\sum_{j_1=0}^{\infty}\ldots
\sum_{j_{s-1}=0}^{\infty}
\sum_{j_s=p+1}^{\infty}\sum_{j_{s+1}=0}^{\infty}\ldots \sum_{j_k=0}^{\infty}
C_{j_k\ldots j_1}^2,
$$

\noindent
where $m_1,\ldots,m_{s-1}>p.$

Denote
$$
C_{j_s\ldots j_1}(\tau)=
\int\limits_t^{\tau}
\phi_{j_s}(t_s)\psi_s(t_s)\ldots \int\limits_t^{t_2}
\phi_{j_1}(t_1)\psi_1(t_1)dt_1\ldots dt_s,
$$
where
$s=1,\ldots,k-1.$

For $s<k$ due to Lemma 3, Dini Theorem (see Appendix)
and
Parseval equality
we obtain
$$
\sum_{j_1=0}^{\infty}\ldots
\sum_{j_{s-1}=0}^{\infty}
\sum_{j_s=p+1}^{\infty}\sum_{j_{s+1}=0}^{\infty}\ldots \sum_{j_k=0}^{\infty}
C_{j_k\ldots j_1}^2=
$$
$$
=
\sum_{j_s=p+1}^{\infty}
\sum_{j_{s-1}=0}^{\infty}\ldots
\sum_{j_{1}=0}^{\infty}
\sum_{j_{s+1}=0}^{\infty}\ldots \sum_{j_k=0}^{\infty}
C_{j_k\ldots j_1}^2=
$$
$$
=
\sum_{j_s=p+1}^{\infty}
\sum_{j_{s-1}=0}^{\infty}\ldots
\sum_{j_{1}=0}^{\infty}
\sum_{j_{s+1}=0}^{\infty}\ldots 
\sum_{j_{k-1}=0}^{\infty}
\int\limits_t^T \psi_k^2(t_k) \left(C_{j_{k-1}\ldots j_1}(t_k)\right)^2 dt_k=
$$
$$
=
\sum_{j_s=p+1}^{\infty}
\sum_{j_{s-1}=0}^{\infty}\ldots
\sum_{j_{1}=0}^{\infty}
\sum_{j_{s+1}=0}^{\infty}\ldots 
\sum_{j_{k-2}=0}^{\infty}
\int\limits_t^T \psi_k^2(t_k) 
\sum_{j_{k-1}=0}^{\infty}\left(C_{j_{k-1}\ldots j_1}(t_k)\right)^2 dt_k=
$$
$$
=
\sum_{j_s=p+1}^{\infty}
\sum_{j_{s-1}=0}^{\infty}\ldots
\sum_{j_{1}=0}^{\infty}
\sum_{j_{s+1}=0}^{\infty}\ldots 
\sum_{j_{k-2}=0}^{\infty}
\int\limits_t^T \psi_k^2(t_k) \int\limits_t^{t_k} \psi_{k-1}^2(\tau) 
\left(C_{j_{k-2}\ldots j_1}(\tau)\right)^2
d\tau dt_k\le
$$

$$
\le M
\sum_{j_s=p+1}^{\infty}
\sum_{j_{s-1}=0}^{\infty}\ldots
\sum_{j_{1}=0}^{\infty}
\sum_{j_{s+1}=0}^{\infty}\ldots 
\sum_{j_{k-2}=0}^{\infty}
\int\limits_t^T  
\left(C_{j_{k-2}\ldots j_1}(\tau)\right)^2 d\tau=
$$
$$
= M
\sum_{j_s=p+1}^{\infty}
\sum_{j_{s-1}=0}^{\infty}\ldots
\sum_{j_{1}=0}^{\infty}
\sum_{j_{s+1}=0}^{\infty}\ldots 
\sum_{j_{k-3}=0}^{\infty}
\int\limits_t^T 
\sum_{j_{k-2}=0}^{\infty}
\left(C_{j_{k-2}\ldots j_1}(\tau)\right)^2
d\tau= 
$$
\newpage
\noindent
$$
= M
\sum_{j_s=p+1}^{\infty}
\sum_{j_{s-1}=0}^{\infty}\ldots
\sum_{j_{1}=0}^{\infty}
\sum_{j_{s+1}=0}^{\infty}\ldots 
\sum_{j_{k-3}=0}^{\infty}
\int\limits_t^T \int\limits_t^{\tau}
\psi_{k-2}^2(\theta)
\left(C_{j_{k-3}\ldots j_1}(\theta)\right)^2
d\theta d\tau\le 
$$
$$
\le M'
\sum_{j_s=p+1}^{\infty}
\sum_{j_{s-1}=0}^{\infty}\ldots
\sum_{j_{1}=0}^{\infty}
\sum_{j_{s+1}=0}^{\infty}\ldots 
\sum_{j_{k-3}=0}^{\infty}
\int\limits_t^T
\left(C_{j_{k-3}\ldots j_1}(\tau)\right)^2
d\tau\le \ldots \le
$$
$$
\le M_k
\sum_{j_s=p+1}^{\infty}
\sum_{j_{s-1}=0}^{\infty}\ldots
\sum_{j_{1}=0}^{\infty}
\int\limits_t^T 
\left(C_{j_{s}\ldots j_1}(\tau)\right)^2 d\tau=
$$
\begin{equation}
\label{d14}
=M_k
\sum_{j_s=p+1}^{\infty}
\sum_{j_{s-1}=0}^{\infty}\ldots
\sum_{j_{2}=0}^{\infty}
\int\limits_t^T  \sum_{j_{1}=0}^{\infty}
\left(C_{j_{s}\ldots j_1}(\tau)\right)^2 d\tau,
\end{equation}

\vspace{3mm}
\noindent
where constants $M,$ $M'$ depend on $T-t$ and
constant $M_k$ depends on $T-t$ and $k.$

Let us explane more precisely how we obtain (\ref{d14}).
For any function $g(s)\in L_2([t,T])$ we have the following
Parseval equality
$$
\sum\limits_{j=0}^{\infty}\left(\int\limits_t^{\tau}
\phi_j(s)g(s)ds\right)^2=
\sum\limits_{j=0}^{\infty}\left(\int\limits_t^T
{\bf 1}_{\{s<\tau\}}\phi_j(s)g(s)ds\right)^2=
$$
\begin{equation}
\label{d15}
=\int\limits_t^T
\left({\bf 1}_{\{s<\tau\}}\right)^2 g^2(s)ds=
\int\limits_t^{\tau}
g^2(s)ds.
\end{equation}

Equality (\ref{d15}) has been applied repeatedly when we obtaining
(\ref{d14}).

Using the replacement of integration order for Riemann integrals, we have
$$
C_{j_s\ldots j_1}(\tau)=
\int\limits_t^{\tau}
\phi_{j_s}(t_s)\psi_s(t_s)\ldots \int\limits_t^{t_2}
\phi_{j_1}(t_1)\psi_1(t_1)dt_1\ldots dt_s=
$$
$$
=\int\limits_t^{\tau}
\phi_{j_1}(t_1)\psi_1(t_1)\int\limits_{t_1}^{\tau}
\phi_{j_2}(t_2)\psi_2(t_2)
\ldots
\int\limits_{t_{s-1}}^{\tau}
\phi_{j_s}(t_s)\psi_s(t_s)dt_s\ldots dt_2dt_1.
$$

\vspace{2mm}

For $l=1,\ldots,s$ we will use the following notation

\newpage
\noindent
$$
{\tilde C}_{j_s\ldots j_l}(\tau,\theta)=
$$
$$
=
\int\limits_{\theta}^{\tau}
\phi_{j_l}(t_l)\psi_l(t_l)\int\limits_{t_l}^{\tau}
\phi_{j_{l+1}}(t_{l+1})\psi_{l+1}(t_{l+1})
\ldots
\int\limits_{t_{s-1}}^{\tau}
\phi_{j_s}(t_s)\psi_s(t_s)dt_s\ldots dt_{l+1}dt_l.
$$

\vspace{2mm}

Using the Parseval equality and Dini Theorem (see
Appendix), from (\ref{d14}) we obtain

\vspace{-2mm}
$$
\sum_{j_1=0}^{\infty}\ldots
\sum_{j_{s-1}=0}^{\infty}
\sum_{j_s=p+1}^{\infty}\sum_{j_{s+1}=0}^{\infty}\ldots \sum_{j_k=0}^{\infty}
C_{j_k\ldots j_1}^2\le
$$
$$
\le
M_k
\sum_{j_s=p+1}^{\infty}
\sum_{j_{s-1}=0}^{\infty}\ldots
\sum_{j_{2}=0}^{\infty}
\int\limits_t^T  \sum_{j_{1}=0}^{\infty}
\left(C_{j_{s}\ldots j_1}(\tau)\right)^2 d\tau=
$$
\begin{equation}
\label{molod1}
=M_k
\sum_{j_s=p+1}^{\infty}
\sum_{j_{s-1}=0}^{\infty}\ldots
\sum_{j_{2}=0}^{\infty}
\int\limits_t^T\int\limits_t^{\tau}\psi_1^2(t_1)  
\left({\tilde C}_{j_{s}\ldots j_2}(\tau,t_1)\right)^2 dt_1d\tau=
\end{equation}
\begin{equation}
\label{molod2}
~~~~~=M_k
\sum_{j_s=p+1}^{\infty}
\sum_{j_{s-1}=0}^{\infty}\ldots
\sum_{j_{3}=0}^{\infty}
\int\limits_t^T\int\limits_t^{\tau}\psi_1^2(t_1)  
\sum_{j_{2}=0}^{\infty}
\left({\tilde C}_{j_{s}\ldots j_2}(\tau,t_1)\right)^2 dt_1d\tau=
\end{equation}
$$
=M_k
\sum_{j_s=p+1}^{\infty}
\sum_{j_{s-1}=0}^{\infty}\ldots
\sum_{j_{3}=0}^{\infty}
\int\limits_t^T\int\limits_t^{\tau}\psi_1^2(t_1)  
\int\limits_{t_1}^{\tau}\psi_2^2(t_2)  
\left({\tilde C}_{j_{s}\ldots j_3}(\tau,t_2)\right)^2 dt_2dt_1d\tau\le
$$
$$
\le M_k
\sum_{j_s=p+1}^{\infty}
\sum_{j_{s-1}=0}^{\infty}\ldots
\sum_{j_{3}=0}^{\infty}
\int\limits_t^T\int\limits_t^{\tau}\psi_1^2(t_1)  
\int\limits_{t}^{\tau}\psi_2^2(t_2)  
\left({\tilde C}_{j_{s}\ldots j_3}(\tau,t_2)\right)^2 dt_2dt_1d\tau\le
$$
$$
\le M^{'}_k
\sum_{j_s=p+1}^{\infty}
\sum_{j_{s-1}=0}^{\infty}\ldots
\sum_{j_{3}=0}^{\infty}
\int\limits_t^T
\int\limits_{t}^{\tau}\psi_2^2(t_2)  
\left({\tilde C}_{j_{s}\ldots j_3}(\tau,t_2)\right)^2 dt_2d\tau
\le \ldots \le
$$
$$
\le M^{''}_k
\sum_{j_s=p+1}^{\infty}
\int\limits_t^T\int\limits_t^{\tau}
\psi_{s-1}^2(t_{s-1})
\left({\tilde C}_{j_{s}}(\tau,t_{s-1})\right)^2 dt_{s-1} d\tau\le
$$
\begin{equation}
\label{la}
\le {\tilde M}_k
\sum_{j_s=p+1}^{\infty}
\int\limits_t^T\int\limits_t^{\tau}
\left(~\int\limits_{u}^{\tau}\phi_{j_s}(\theta)
\psi_s(\theta)d\theta\right)^2 du d\tau,
\end{equation}

\newpage
\noindent
where constants $M_k^{'},$ $M_k^{''},$
and ${\tilde M}_k$ depend on $k$ and $T-t.$

Let us explane more precisely how we obtain (\ref{la}).
For any function $g(s)\in L_2([t,T])$ we have the following
Parseval equality
$$
\sum\limits_{j=0}^{\infty}\left(\int\limits_{\theta}^{\tau}
\phi_j(s)g(s)ds\right)^2=
\sum\limits_{j=0}^{\infty}\left(\int\limits_t^T
{\bf 1}_{\{\theta<s<\tau\}}\phi_j(s)g(s)ds\right)^2=
$$
\begin{equation}
\label{d22}
=\int\limits_t^T
\left({\bf 1}_{\{\theta<s<\tau\}}\right)^2 g^2(s)ds=
\int\limits_{\theta}^{\tau}
g^2(s)ds.
\end{equation}

Equality (\ref{d22}) has been applied repeatedly when we obtain
(\ref{la}).

Let us explane more precisely the passing from (\ref{molod1})
to (\ref{molod2}) (the same steps have been used when we 
derived (\ref{la})). 

We have
$$
\int\limits_t^T\int\limits_t^{\tau}\psi_1^2(t_1)  
\sum_{j_{2}=0}^{\infty}
\left({\tilde C}_{j_{s}\ldots j_2}(\tau,t_1)\right)^2 dt_1d\tau -
\sum_{j_{2}=0}^{n}\int\limits_t^T\int\limits_t^{\tau}\psi_1^2(t_1)  
\left({\tilde C}_{j_{s}\ldots j_2}(\tau,t_1)\right)^2 dt_1d\tau =
$$
$$
=\int\limits_t^T\int\limits_t^{\tau}\psi_1^2(t_1)  
\sum_{j_{2}=n+1}^{\infty}
\left({\tilde C}_{j_{s}\ldots j_2}(\tau,t_1)\right)^2 dt_1d\tau =
$$
\begin{equation}
\label{molod3}
=\lim\limits_{N\to\infty}
\sum\limits_{j=0}^{N-1}\int\limits_t^{\tau_j}\psi_1^2(t_1)  
\sum_{j_{2}=n+1}^{\infty}
\left({\tilde C}_{j_{s}\ldots j_2}(\tau_j,t_1)\right)^2 dt_1 \Delta\tau_j,
\end{equation}

\noindent
where $\{\tau_j\}_{j=0}^{N}$ is the partition of the 
interval $[t, T],$ which satisfies the condition (\ref{1111}).

Since the non-decreasing functional sequence $u_n(\tau_j,t_1)$ and its
limit function $u(\tau_j,t_1)$ are continuous on the
interval $[t,\tau_j]\subseteq [t, T]$ with respect to $t_1$,
where
$$
u_n(\tau_j,t_1)=
\sum_{j_{2}=0}^{n}
\left({\tilde C}_{j_{s}\ldots j_2}(\tau_j,t_1)\right)^2,
$$
$$
u(\tau_j,t_1)=
\sum_{j_{2}=0}^{\infty}
\left({\tilde C}_{j_{s}\ldots j_2}(\tau_j,t_1)\right)^2=
\int\limits_{t_1}^{\tau_j}
\psi_2^2(t_2)
\left({\tilde C}_{j_{s}\ldots j_3}(\tau_j,t_2)\right)^2 dt_2,
$$

\noindent 
then by Dini Theorem we have the uniform convergence
of $u_n(\tau_j,t_1)$ to $u(\tau_j,t_1)$ at the interval $[t,\tau_j]\subseteq
[t, T]$
with respect to $t_1.$ As a result, we obtain
\begin{equation}
\label{molod4}
\sum_{j_{2}=n+1}^{\infty}
\left({\tilde C}_{j_{s}\ldots j_2}(\tau_j,t_1)\right)^2<\varepsilon,\ \ \ 
t_1\in [t,\tau_j]
\end{equation}

\noindent
for $n>N(\varepsilon)$ ($N(\varepsilon)$ exists
for any $\varepsilon>0$ and it does not depend on $t_1$).

From (\ref{molod3}) and (\ref{molod4}) we obtain
$$
\lim\limits_{N\to\infty}
\sum\limits_{j=0}^{N-1}\int\limits_t^{\tau_j}\psi_1^2(t_1)  
\sum_{j_{2}=n+1}^{\infty}
\left({\tilde C}_{j_{s}\ldots j_2}(\tau_j,t_1)\right)^2 dt_1 \Delta\tau_j
\le
$$
$$
\le
\varepsilon 
\lim\limits_{N\to\infty}
\sum\limits_{j=0}^{N-1}\int\limits_t^{\tau_j}\psi_1^2(t_1)  
dt_1 \Delta\tau_j= 
$$
\begin{equation}
\label{molod6}
=\varepsilon \int\limits_t^T
\int\limits_t^{\tau}\psi_1^2(t_1)  
dt_1 d\tau.
\end{equation}

From (\ref{molod6}) we get
$$
\lim\limits_{n\to\infty}\int\limits_t^T\int\limits_t^{\tau}\psi_1^2(t_1)  
\sum_{j_{2}=n+1}^{\infty}
\left({\tilde C}_{j_{s}\ldots j_2}(\tau,t_1)\right)^2 dt_1d\tau = 0.
$$

This fact completes the proof of passing 
from (\ref{molod1})
to (\ref{molod2}).

Let us estimate the integral 
\begin{equation}
\label{st1}
\int\limits_{u}^{\tau}\phi_{j_s}(\theta)
\psi_s(\theta)d\theta
\end{equation}
from (\ref{la}) for the cases when $\{\phi_j(s)\}_{j=0}^{\infty}$
is a complete orthonormal system of Legendre polynomials or
trigonometric functions in the space $L_2([t,T])$.

Note that the estimates for the integral
\begin{equation}
\label{st2}
\int\limits_{t}^{\tau}\phi_{j}(\theta)\psi(\theta)d\theta,\ \ \ j\ge p+1
\end{equation}
have been obtained in \cite{Kuznetsov6}, \cite{Kuznetsov7},
\cite{Kuznetsov8}, \cite{Kuznetsov16}. Here
$\psi(\theta)$ is a continuously
differentiable function on the interval $[t, T]$,

Let us estimate the integral (\ref{st1}) using the approach from
\cite{Kuznetsov6}, \cite{Kuznetsov7},
\cite{Kuznetsov8}, \cite{Kuznetsov16}.

First consider the case of Legendre polynomials.
Then $\phi_j(\theta)$ looks as follows

\vspace{-2mm}
$$
\phi_j(\theta)=\sqrt{\frac{2j+1}{T-t}}P_j\left(\left(
\theta-\frac{T+t}{2}\right)\frac{2}{T-t}\right),\ \ \ j\ge 0,
$$

\vspace{2mm}
\noindent
where $P_j(x)$ $(j=0, 1, 2\ldots)$ is the 
Legendre polynomial.

Further, we have 
$$
\int\limits_v^x\phi_{j}(\theta)\psi(\theta)d\theta=
\frac{\sqrt{T-t}\sqrt{2j+1}}{2}
\int\limits_{z(v)}^{z(x)}P_{j}(y)
\psi(u(y))dy=
$$
$$
=\frac{\sqrt{T-t}}{2\sqrt{2j+1}}\Biggl((P_{j+1}(z(x))-
P_{j-1}(z(x)))\psi(x)-
$$
$$
-
(P_{j+1}(z(v))-
P_{j-1}(z(v)))\psi(v)-
\Biggr.
$$
\begin{equation}
\label{6000}
\Biggl.-
\frac{T-t}{2}
\int\limits_{z(v)}^{z(x)}((P_{j+1}(y)-P_{j-1}(y))
{\psi}'(u(y))dy\Biggr),
\end{equation}

\noindent
where $x, v\in (t, T),$ $j\ge p+1,$ and
$u(y)$, $z(x)$ are defined by the following relations
$$
u(y)=\frac{T-t}{2}y+\frac{T+t}{2},\ \ \
z(x)=\left(x-\frac{T+t}{2}\right)\frac{2}{T-t},
$$

\noindent
${\psi}'$ is a derivative of the function $\psi(\theta)$
with respect to the variable $u(y).$

Note that in (\ref{6000}) we used the following well-known property
of the Legendre polynomials \cite{Suetin}
$$
\frac{dP_{j+1}}{dx}(x)-\frac{dP_{j-1}}{dx}(x)=(2j+1)P_j(x),\ \ \ 
j=1, 2,\ldots
$$

\vspace{1mm}

From (\ref{6000}) and the well-known estimate for the Legendre
polynomials \cite{Gob}
$$
|P_j(y)| <\frac{K}{\sqrt{j+1}(1-y^2)^{1/4}},\ \ \ 
y\in (-1, 1),\ \ \ j\in \mathbb{N},
$$

\vspace{1mm}
\noindent
where constant $K$ does not depend on $y$ and $j$, it follows that
\begin{equation}
\label{101oh}
\left|
\int\limits_v^x\phi_{j}(\theta)\psi(\theta)d\theta
\right| <
\frac{C}{j}\Biggl(\frac{1}{(1-(z(x))^2)^{1/4}}+
\frac{1}{(1-(z(v))^2)^{1/4}}+C_1\Biggr),
\end{equation}

\noindent
where $z(x), z(v)\in (-1, 1),$ $x, v\in (t, T)$ and
constants $C, C_1$ does not depend on $j$.

From (\ref{101oh}) we obtain
\begin{equation}
\label{102oh}
\left(
\int\limits_v^x\phi_{j}(\theta)\psi(\theta)d\theta
\right)^2 
<
\frac{C_2}{j^2}\Biggl(\frac{1}{(1-(z(x))^2)^{1/2}}+
\frac{1}{(1-(z(v))^2)^{1/2}}+C_3\Biggr),
\end{equation}
where constants $C_2, C_3$ does not depend on $j$.

Let us apply (\ref{102oh}) for the estimate of the right-hand side
of (\ref{la}). We have
$$
\int\limits_t^T\int\limits_t^{\tau}
\left(~\int\limits_{u}^{\tau}\phi_{j_s}(\theta)
\psi_s(\theta)d\theta\right)^2 du d\tau\le
$$
$$
\le \frac{K_1}{j_s^2}
\left(
\int\limits_{-1}^1
\frac{dy}{\left(1-y^2\right)^{1/2}}+
\int\limits_{-1}^1\int\limits_{-1}^x
\frac{dy}{\left(1-y^2\right)^{1/2}}dx + K_2\right)\le
$$

\begin{equation}
\label{103}
\le\frac{K_3}{j_s^2},
\end{equation}

\vspace{2mm}
\noindent
where constants $K_1, K_2, K_3$ are independent of $j_s.$

Now consider the trigonometric case.
The complete orthonormal system of trigonometric functions
in the space $L_2([t, T])$ has the following form
\begin{equation}
\label{trig11oh}
\phi_j(\theta)=\frac{1}{\sqrt{T-t}}
\left\{
\begin{matrix}
1,\ & j=0\cr\cr
\sqrt{2}{\rm sin} \left(2\pi r(\theta-t)/(T-t)\right),\ & j=2r-1\cr\cr
\sqrt{2}{\rm cos} \left(2\pi r(\theta-t)/(T-t)\right),\ & j=2r
\end{matrix}
,\right.
\end{equation}

\noindent
where $r=1, 2,\ldots $

Using the system of functions 
(\ref{trig11oh}), we have
$$
\int\limits_v^x\phi_{2r-1}(\theta)\psi(\theta)d\theta=
\sqrt{\frac{2}{T-t}}\int\limits_v^x
{\rm sin} \frac{2\pi r(\theta-t)}{T-t}\psi(\theta)d\theta=
$$
$$
=-\sqrt{\frac{T-t}{2}}\frac{1}{\pi r}\Biggl(
\psi(x){\rm cos}\frac{2\pi r(x-t)}{T-t}-
\psi(v){\rm cos}\frac{2\pi r(v-t)}{T-t}-\Biggr.
$$
\begin{equation}
\label{201}
\Biggl.-
\int\limits_v^x
{\rm cos} \frac{2\pi r(\theta-t)}{T-t}\psi'(\theta)d\theta\Biggr),
\end{equation}

$$
\int\limits_v^x\phi_{2r}(\theta)\psi(\theta)d\theta=
\sqrt{\frac{2}{T-t}}\int\limits_v^x
{\rm cos} \frac{2\pi r(\theta-t)}{T-t}\psi(\theta)d\theta=
$$
$$
=\sqrt{\frac{T-t}{2}}\frac{1}{\pi r}\Biggl(
\psi(x){\rm sin}\frac{2\pi r(x-t)}{T-t}-
\psi(v){\rm sin}\frac{2\pi r(v-t)}{T-t}-\Biggr.
$$
\begin{equation}
\label{202}
\Biggl.-
\int\limits_v^x
{\rm sin} \frac{2\pi r(\theta-t)}{T-t}\psi'(\theta)d\theta\Biggr),
\end{equation}

\vspace{2mm}
\noindent
where $\psi'(\theta)$ is a derivative of the function $\psi(\theta)$
with respect to the variable $\theta.$

Combining (\ref{201}) and (\ref{202}), we obtain for the
trigonometric case
\begin{equation}
\label{203}
\left(
\int\limits_v^x\phi_{j}(\theta)\psi(\theta)d\theta
\right)^2 \le 
\frac{C_4}{j^2},
\end{equation}

\noindent
where constant $C_4$ is independent of $j.$

From (\ref{203}) we finally have
\begin{equation}
\label{103x}
\int\limits_t^T\int\limits_t^{\tau}
\left(~\int\limits_{u}^{\tau}\phi_{j_s}(\theta)
\psi_s(\theta)d\theta\right)^2 du d\tau
\le \frac{K_4}{j_s^2},
\end{equation}

\noindent
where constant $K_4$ is independent of $j_s.$

Combining (\ref{la}), (\ref{103}) and (\ref{103x}), we obtain
$$
\sum_{j_1=0}^{\infty}\ldots
\sum_{j_{s-1}=0}^{\infty}
\sum_{j_s=p+1}^{\infty}\sum_{j_{s+1}=0}^{\infty}\ldots \sum_{j_k=0}^{\infty}
C_{j_k\ldots j_1}^2\le
$$
\begin{equation}
\label{fffoh}
\le L_k
\sum_{j_s=p+1}^{\infty}\frac{1}{j_s^2} \le 
L_k \int\limits_{p}^{\infty}\frac{dx}{x^2}=\frac{L_k}{p},
\end{equation}

\noindent
where constant $L_k$ depends on $k$ and $T-t.$

Obviously, the case $s=k$ can be considered absolutely analogously to the
case $s<k$. Then from (\ref{aaap}) and (\ref{fffoh})
we obtain
\begin{equation}
\label{ddd1}
\int\limits_{[t,T]^k}K^2(t_1,\ldots,t_k)dt_1\ldots dt_k-
\sum_{j_1=0}^{p}\ldots \sum_{j_k=0}^{p}
C_{j_k\ldots j_1}^2\le \frac{G_k}{p},
\end{equation}
where constant $G_k$ depends on $k$ and $T-t.$

For the further consideration we will use estimate (\ref{dima2ye100}).
Using (\ref{ddd1}) and the estimate (\ref{dima2ye100})
for the case $p_1=\ldots=p_k=p$ and $n=2$, 
we obtain
$$
{\sf M}\left\{\biggl(J[\psi^{(k)}]_{T,t}-
J[\psi^{(k)}]_{T,t}^{p,\ldots,p}\biggr)^{4}\right\}\le
$$
$$
\le C_{2,k}
\left(
\int\limits_{[t,T]^k}
K^2(t_1,\ldots,t_k)
dt_1\ldots dt_k -\sum_{j_1=0}^{p}\ldots
\sum_{j_k=0}^{p}C^2_{j_k\ldots j_1}
\right)^2\le 
$$
\begin{equation}
\label{fff5}
\le\frac{H_{2,k}}{p^2},
\end{equation}

\vspace{2mm}
\noindent
where 
$$
C_{n,k}=(k!)^{n} (2n-1)^{nk}
$$

\vspace{2mm}
\noindent
and $H_{2,k}=G_k^2{C}_{2,k}.$

Let us consider Lemma 1 (see Appendix) with
$$
\xi_p=\biggl|J[\psi^{(k)}]_{T,t}-
J[\psi^{(k)}]_{T,t}^{p,\ldots,p}\biggr|\ \ \ \hbox{and}\ \ \ 
\alpha=4.
$$

Then from (\ref{fff5}) we get
\begin{equation}
\label{qqq1oh}
\sum\limits_{p=1}^{\infty}
{\sf M}\left\{\biggl(J[\psi^{(k)}]_{T,t}-
J[\psi^{(k)}]_{T,t}^{p,\ldots,p}\biggr)^{4}\right\}\le
H_{2,k}\sum\limits_{p=1}^{\infty}\frac{1}{p^2}<\infty.
\end{equation}

\vspace{2mm}

Using Lemma 1 (see Appendix) and the estimate (\ref{qqq1oh}), we obtain

\vspace{-1mm}
$$
J[\psi^{(k)}]_{T,t}^{p,\ldots,p}\ \to \ J[\psi^{(k)}]_{T,t}\ \ \ 
\hbox{if}\ \ \ p\to \infty\ \ \ \hbox{w.\ p.\ 1},
$$

\noindent
where (see Theorem 1)
$$
J[\psi^{(k)}]_{T,t}^{p,\ldots,p}=
\sum_{j_1,\ldots,j_k=0}^{p}
C_{j_k\ldots j_1}\Biggl(
\prod_{l=1}^k\zeta_{j_l}^{(i_l)}\ -
\Biggr.
$$
\begin{equation}
\label{kk0}
-\ \Biggl.
\hbox{\vtop{\offinterlineskip\halign{
\hfil#\hfil\cr
{\rm l.i.m.}\cr
$\stackrel{}{{}_{N\to \infty}}$\cr
}} }\sum_{(l_1,\ldots,l_k)\in {\rm G}_k}
\phi_{j_{1}}(\tau_{l_1})
\Delta{\bf w}_{\tau_{l_1}}^{(i_1)}\ldots
\phi_{j_{k}}(\tau_{l_k})
\Delta{\bf w}_{\tau_{l_k}}^{(i_k)}\Biggr)
\end{equation}
or (see Theorem 2)
$$
J[\psi^{(k)}]_{T,t}^{p,\ldots,p}=
\sum\limits_{j_1,\ldots,j_k=0}^{p}
C_{j_k\ldots j_1}\Biggl(
\prod_{l=1}^k\zeta_{j_l}^{(i_l)}+\sum\limits_{r=1}^{[k/2]}
(-1)^r \times
\Biggr.
$$
\begin{equation}
\label{kkohh}
\times
\sum_{\stackrel{(\{\{g_1, g_2\}, \ldots, 
\{g_{2r-1}, g_{2r}\}\}, \{q_1, \ldots, q_{k-2r}\})}
{{}_{\{g_1, g_2, \ldots, 
g_{2r-1}, g_{2r}, q_1, \ldots, q_{k-2r}\}=\{1, 2, \ldots, k\}}}}
\prod\limits_{s=1}^r
{\bf 1}_{\{i_{g_{{}_{2s-1}}}=~i_{g_{{}_{2s}}}\ne 0\}}
\Biggl.{\bf 1}_{\{j_{g_{{}_{2s-1}}}=~j_{g_{{}_{2s}}}\}}
\prod_{l=1}^{k-2r}\zeta_{j_{q_l}}^{(i_{q_l})}\Biggr),
\end{equation}

\vspace{3mm}
\noindent
where $i_1,\ldots,i_k=1,\ldots,m$ in (\ref{kk0}) and (\ref{kkohh}).
The proof of Theorem 6 is completed.

\section{Appendix}

~~~~{\bf Lemma 1} \cite{Shir}. {\it If for the sequence of random variables
$\xi_p$ and for some
$\alpha>0$ the number series 
$$
\sum\limits_{p=1}^{\infty}{\sf M}\left\{\left|\xi_p\right|^{\alpha}\right\}
$$
converges, then the sequence $\xi_p$ converges to zero w.~p.~{\rm 1}.}

{\bf Lemma 2.}\ {\it The following equalities are fulfilled
$$
\sum_{j_1,\ldots,j_k=0}^{\infty}
C_{j_k\ldots j_1}^2=
\sum_{j_1=0}^{\infty}\ldots
\sum_{j_k=0}^{\infty}
C_{j_k\ldots j_1}^2=
$$
\begin{equation}
\label{lem1}
=\sum_{j_k=0}^{\infty}\ldots
\sum_{j_1=0}^{\infty}
C_{j_k\ldots j_1}^2=
\sum_{j_{q_1}=0}^{\infty}\ldots
\sum_{j_{q_k}=0}^{\infty}
C_{j_k\ldots j_1}^2
\end{equation}

\vspace{1mm}
\noindent
for any permutation $(q_1,\ldots,q_k)$ such that
$\{q_1,\ldots,q_k\}=\{1,\ldots,k\},$
where $C_{j_k\ldots j_1}$ is defined by {\rm (\ref{rrr1})}.
}

{\bf Proof.} Let us remind the well-known 
fact from the mathematical analysis,
which is connected to existence
of iterated limits.

{\bf Proposition 1} \cite{IP}.
{\it Let $\bigl\{x_{n,m}\bigr\}_{n,m=1}^{\infty}$
be a double sequence and let there exists the limit
$$
\lim\limits_{n,m\to\infty}x_{n,m}=a<\infty.
$$

Moreover, let there exist the limits
$$
\lim\limits_{n\to\infty}x_{n,m}<\infty\ \ \ \hbox{for any}\ \ m,\ \ \ \ \
\lim\limits_{m\to\infty}x_{n,m}<\infty\ \ \ \hbox{for any}\ \ n.
$$
Then there exist the iterated limits
$$
\lim\limits_{n\to\infty}\lim\limits_{m\to\infty}x_{n,m},\ \ \
\lim\limits_{m\to\infty}\lim\limits_{n\to\infty}x_{n,m}
$$
and moreover,
$$
\lim\limits_{n\to\infty}\lim\limits_{m\to\infty}x_{n,m}=
\lim\limits_{m\to\infty}\lim\limits_{n\to\infty}x_{n,m}=a.
$$
}

Let us consider the value
\begin{equation}
\label{21}
\sum_{j_{q_l}=0}^{p}\ldots
\sum_{j_{q_k}=0}^{p}
C_{j_k\ldots j_1}^2
\end{equation}
for any permutation $(q_l,\ldots,q_k)$, where $l=1,2,\ldots,k$,
$\{q_1,\ldots,q_k\}=\{1,\ldots,k\}.$

Obviously, (\ref{21}) 
is the non-decreasing sequence with respect to $p$.
Moreover,
$$
\sum_{j_{q_l}=0}^{p}\ldots
\sum_{j_{q_k}=0}^{p}
C_{j_k\ldots j_1}^2\le 
\sum_{j_{q_1}=0}^{p}\sum_{j_{q_2}=0}^{p}\ldots
\sum_{j_{q_k}=0}^{p}
C_{j_k\ldots j_1}^2\le 
$$
$$
\le
\sum_{j_1,\ldots,j_k=0}^{\infty}
C_{j_k\ldots j_1}^2<\infty.
$$

Then the following limit
$$
\lim\limits_{p\to\infty}\sum\limits_{j_{q_l}=0}^p \ldots 
\sum\limits_{j_{q_k}=0}^{p}
C_{j_k\ldots j_1}^2=
\sum_{j_{q_l},\ldots,j_{q_k}=0}^{\infty}
C_{j_k\ldots j_1}^2
$$
exists.

Let $p_l,\ldots,p_k$ simultaneously tend to infinity.
Then $g, r\to \infty$, where $g=\min\{p_l,\ldots,p_k\}$ and
$r=\max\{p_l,\ldots,p_k\}$. Moreover,
$$
\sum_{j_{q_l}=0}^{g}\ldots
\sum_{j_{q_k}=0}^{g}
C_{j_k\ldots j_1}^2\le 
\sum_{j_{q_l}=0}^{p_l}\ldots
\sum_{j_{q_k}=0}^{p_k}
C_{j_k\ldots j_1}^2\le
\sum_{j_{q_l}=0}^{r}\ldots
\sum_{j_{q_k}=0}^{r}
C_{j_k\ldots j_1}^2.
$$

This means that the existence of the limit 
\begin{equation}
\label{1c1c}
\lim\limits_{p\to\infty}\sum_{j_{q_l}=0}^{p}\ldots
\sum_{j_{q_k}=0}^{p}
C_{j_k\ldots j_1}^2
\end{equation}

\noindent
implies the existence of the limit 
\begin{equation}
\label{1d1d}
\lim\limits_{p_l,\ldots,p_k\to\infty}\sum_{j_{q_l}=0}^{p_l}\ldots
\sum_{j_{q_k}=0}^{p_k}
C_{j_k\ldots j_1}^2
\end{equation}
and equality of the limits (\ref{1c1c}) and (\ref{1d1d}).
 
Consequently,
$$
\lim\limits_{p,q\to\infty}\sum_{j_{q_l}=0}^{q}\sum_{j_{q_{l+1}}=0}^{p}\ldots
\sum_{j_{q_k}=0}^{p}
C_{j_k\ldots j_1}^2=
\lim\limits_{p\to\infty}\sum_{j_{q_l}=0}^{p}\ldots
\sum_{j_{q_k}=0}^{p}
C_{j_k\ldots j_1}^2=
$$
\begin{equation}
\label{1h1h}
=\lim\limits_{p_l,\ldots,p_k\to\infty}\sum_{j_{q_l}=0}^{p_l}\ldots
\sum_{j_{q_k}=0}^{p_k}
C_{j_k\ldots j_1}^2.
\end{equation}

Since the limit
$$
\sum_{j_1,\ldots,j_k=0}^{\infty}
C_{j_k\ldots j_1}^2
$$

\noindent
exists (see the Parseval equality (\ref{par1})), then from Proposition 1
we have
$$
\sum_{j_{q_1}=0}^{\infty}\sum_{j_{q_2},\ldots,j_{q_k}=0}^{\infty}
C_{j_k\ldots j_1}^2=
\lim\limits_{q\to\infty}
\lim\limits_{p\to\infty}
\sum_{j_{q_1}=0}^{q}\sum_{j_{q_2}=0}^p \ldots \sum_{j_{q_k}=0}^{p}
C_{j_k\ldots j_1}^2=
$$
\begin{equation}
\label{1b1b}
=\lim\limits_{q,p\to\infty}
\sum_{j_{q_1}=0}^{q}\sum_{j_{q_2}=0}^p \ldots \sum_{j_{q_k}=0}^{p}
C_{j_k\ldots j_1}^2=
\sum_{j_1,\ldots,j_k=0}^{\infty}
C_{j_k\ldots j_1}^2.
\end{equation}

\noindent
\par
Using (\ref{1h1h}) and Proposition 1, we have
$$
\sum_{j_{q_2}=0}^{\infty}\sum_{j_{q_3},\ldots,j_{q_k}=0}^{\infty}
C_{j_k\ldots j_1}^2=
\lim\limits_{q\to\infty}
\lim\limits_{p\to\infty}
\sum_{j_{q_2}=0}^{q}\sum_{j_{q_3}=0}^p \ldots \sum_{j_{q_k}=0}^{p}
C_{j_k\ldots j_1}^2=
$$
\begin{equation}
\label{1a1a}
=\lim\limits_{q,p\to\infty}
\sum_{j_{q_2}=0}^{q}\sum_{j_{q_3}=0}^p \ldots \sum_{j_{q_k}=0}^{p}
C_{j_k\ldots j_1}^2=
\sum_{j_{q_2},\ldots,j_{q_k}=0}^{\infty}
C_{j_k\ldots j_1}^2.
\end{equation}

Combining (\ref{1a1a}) and (\ref{1b1b}), we obtain
$$
\sum_{j_{q_1}=0}^{\infty}\sum_{j_{q_2}=0}^{\infty}
\sum_{j_{q_3},\ldots,j_{q_k}=0}^{\infty}
C_{j_k\ldots j_1}^2=
\sum_{j_{1},\ldots,j_{k}=0}^{\infty}
C_{j_k\ldots j_1}^2.
$$

Repeating the above steps, we complete the proof of Lemma
2.

{\bf Lemma 3.} {\it The following equality takes place
\begin{equation}
\label{d11oh}
\sum_{j_1=0}^{\infty}\ldots
\sum_{j_{s-1}=0}^{\infty}
\sum_{j_s=p+1}^{\infty}\sum_{j_{s+1}=0}^{\infty}\ldots \sum_{j_k=0}^{\infty}
C_{j_k\ldots j_1}^2=
$$
$$
=
\sum_{j_s=p+1}^{\infty}\sum_{j_{s-1}=0}^{\infty}\ldots
\sum_{j_{1}=0}^{\infty}
\sum_{j_{s+1}=0}^{\infty}\ldots \sum_{j_k=0}^{\infty}
C_{j_k\ldots j_1}^2,
\end{equation}

\vspace{2mm}
\noindent
where $s=1,\ldots,k$ and 
$C_{j_k\ldots j_1}$ is defined by {\rm (\ref{rrr1})}.}

{\bf Proof.} Applying the arguments that we used in the proof of 
Lemma 2, we obtain
$$
\lim\limits_{n\to\infty}
\sum_{j_1=0}^{n}\ldots
\sum_{j_{s-1}=0}^{n}
\sum_{j_s=0}^{p}\sum_{j_{s+1}=0}^{n}\ldots \sum_{j_k=0}^{n}
C_{j_k\ldots j_1}^2=
$$
\begin{equation}
\label{ura0}
=\sum_{j_s=0}^{p}\ \sum_{j_{1},\ldots, j_{s-1}, 
j_{s+1},\ldots,j_k=0}^{\infty}
C_{j_k\ldots j_1}^2
=\sum_{j_s=0}^{p}\sum_{j_{q_1}=0}^{\infty}\ldots
\sum_{j_{q_{k-1}}=0}^{\infty}
C_{j_k\ldots j_1}^2
\end{equation}

\vspace{2mm}
\noindent
for any permutation $(q_1,\ldots,q_{k-1})$ such that
$\{q_1,\ldots,q_{k-1}\}=\{1,\ldots,s-1,s+1,\ldots,k\}$,
where $p$ is a fixed natural number.

Obviously, we have
$$
\sum_{j_s=0}^{p}\sum_{j_{q_1}=0}^{\infty}\ldots
\sum_{j_{q_{k-1}}=0}^{\infty}
C_{j_k\ldots j_1}^2=
\sum_{j_{q_1}=0}^{\infty}\ldots \sum_{j_s=0}^{p} \ldots
\sum_{j_{q_{k-1}}=0}^{\infty}C_{j_k\ldots j_1}^2= \ldots =
$$
\begin{equation}
\label{ura1}
=
\sum_{j_{q_1}=0}^{\infty}\ldots 
\sum_{j_{q_{k-1}}=0}^{\infty}
\sum_{j_s=0}^{p}
C_{j_k\ldots j_1}^2.
\end{equation}

\vspace{2mm}

Using (\ref{ura0}), (\ref{ura1}) and Lemma 2, we obtain

\vspace{-1mm}
$$
\sum_{j_1=0}^{\infty}\ldots
\sum_{j_{s-1}=0}^{\infty}
\sum_{j_s=p+1}^{\infty}\sum_{j_{s+1}=0}^{\infty}\ldots \sum_{j_k=0}^{\infty}
C_{j_k\ldots j_1}^2=
$$

\vspace{-2mm}
$$
=
\sum_{j_1=0}^{\infty}\ldots
\sum_{j_{s-1}=0}^{\infty}
\sum_{j_s=0}^{\infty}\sum_{j_{s+1}=0}^{\infty}\ldots \sum_{j_k=0}^{\infty}
C_{j_k\ldots j_1}^2
-\sum_{j_1=0}^{\infty}\ldots
\sum_{j_{s-1}=0}^{\infty}
\sum_{j_s=0}^{p}\sum_{j_{s+1}=0}^{\infty}\ldots \sum_{j_k=0}^{\infty}
C_{j_k\ldots j_1}^2=
$$
$$
=
\sum_{j_s=0}^{\infty}
\sum_{j_{s-1}=0}^{\infty}\ldots
\sum_{j_1=0}^{\infty}\sum_{j_{s+1}=0}^{\infty}\ldots \sum_{j_k=0}^{\infty}
C_{j_k\ldots j_1}^2-
\sum_{j_s=0}^{p}
\sum_{j_{s-1}=0}^{\infty}\ldots
\sum_{j_1=0}^{\infty}\sum_{j_{s+1}=0}^{\infty}\ldots \sum_{j_k=0}^{\infty}
C_{j_k\ldots j_1}^2=
$$
$$
=
\sum_{j_s=p+1}^{\infty}
\sum_{j_{s-1}=0}^{\infty}\ldots
\sum_{j_1=0}^{\infty}\sum_{j_{s+1}=0}^{\infty}\ldots \sum_{j_k=0}^{\infty}
C_{j_k\ldots j_1}^2.
$$

\vspace{2mm}

The equality (\ref{d11oh}) is proved.

{\bf Theorem (Dini)} \cite{IP2}. {\it 
Let the functional sequence $u_n(x)$ 
be non-decreasing at each point of the interval $[a, b]$.
In addition, all the functions $u_n(x)$
of this sequence and the limit function $u(x)$ are continuous on the interval
$[a, b].$ Then the convergence $u_n(x)$ to 
$u(x)$ is uniform on the interval $[a,b].$}

\renewcommand{\refname}{\rm
{\bf Bibliography}}

\end{document}